\documentstyle[11pt]{article}

       \font\tenmsb=msbm10
       \font\sevenmsb=msbm7
       \font\fivemsb=msbm5
       \catcode`\@=11
       \ifx\amstexloaded@\relax\catcode`\@=\active
       \endinput\else\let\amstexloaded@\relax\fi
       \def\spaces@{\space\space\space\space\space}
       \def\spaces@@{\spaces@\spaces@\spaces@\spaces@\spaces@}
       \def\space@.{\futurelet\space@\relax}
       \space@. %
       \def\Err@#1{\errhelp\defaulthelp@\errmessage{AmS-TeX error: #1}}
       \def\relaxnext@{\let\next\relax}
       \def\accentfam@{7}
       
       \def\noaccents@{\def\accentfam@{0}}
       \def\Cal{\relaxnext@\ifmmode\let\next\Cal@\else
       \def\next{\Err@{Use \string\Cal\space only in math mode}}\fi\next}
       \def\Cal@#1{{\Cal@@{#1}}}
       \def\Cal@@#1{\noaccents@\fam\tw@#1}

       \def\Bbb{\relaxnext@\ifmmode\let\next\Bbb@\else
       \def\next{\Err@{Use \string\Bbb\space only in math mode}}\fi\next}
       
       \def\Bbb@#1{{\Bbb@@{#1}}}
       \def\Bbb@@#1{\noaccents@\fam\msbfam#1}
       \newfam\msbfam
       \textfont\msbfam=\tenmsb
       \scriptfont\msbfam=\sevenmsb
       \scriptscriptfont\msbfam=\fivemsb

\def\N{{\Bbb N}}
\def\Z{{\Bbb Z}}
\def\R{{\Bbb R}}
\def\T{{\Bbb T}}

\newtheorem{Theorem}{Theorem}
\newtheorem{Lemma}{Lemma}[section]
\newtheorem{Proposition}{Proposition}

\@addtoreset{equation}{section}

\newcommand{\beq}{\begin{equation} }
\newcommand{\eeq}{\end{equation} }

\makeatletter
\newcommand{\rmnum}[1]{\romannumeral #1}
\newcommand{\Rmnum}[1]{\expandafter\@slowromancap\romannumeral #1@}
\makeatother

\begin{document}

\setlength{\columnsep}{5pt}

\title{An infinite dimensional KAM theorem with application to two dimensional completely resonant beam equation \thanks{This work is partially supported
by NSFC grant 11271180.}}
\author{\\Jiansheng Geng, Shidi Zhou\\
 {\footnotesize Department of Mathematics,Nanjing University, Nanjing 210093, P.R.China}\\
{\footnotesize Email: jgeng@nju.edu.cn; mathsdzhou@126.com}}


\date{}
\maketitle

\begin{abstract}
 In this paper we consider the completely resonant beam equation on $\T^2$ with cubic nonlinearity on a subspace of $L^2 (\T^2)$ which will be explained later. We establish an abstract infinite dimensional KAM theorem and apply it to the completely resonant beam equation. We prove the existence of a class of Whitney smooth small amplitude quasi-periodic solutions corresponding to finite dimensional tori.
\end{abstract}



\section{\textbf{Introduction}}
$\\$
\indent  In this paper we consider the two dimensional completely resonant beam equation with cubic nonlinearity on a subspace $\mathcal U$ of $L^2 (\T^2)$:

\begin{eqnarray}
u_{tt} + \Delta^2 _x u + u^3 = 0 \quad u=u(t,x),t\in \R, x\in \T^2
\end{eqnarray}
\noindent Here $t$ is time and $x$ is the spatial variable. The subspace $\mathcal U$ is defined by
\begin{eqnarray}
{\mathcal U} =
\{
u=\sum_{n\in \Z^2 _{odd}}u_n \phi_n,\quad \phi_n(x) = e^{i\left\langle n,x \right\rangle}
\}
\end{eqnarray}
\noindent where the integer set $\Z^2 _{odd}$ is defined as
\begin{eqnarray}
\Z^2 _{odd} = \{ n=(n_1, n_2) : n_1 \in 2\Z-1, n_2 \in 2\Z  \}
\end{eqnarray}
\noindent This idea comes from the work by M.Procesi [29] and we
will explain it later in section 2. The solution of "real"
completely resonant beam equation (not on $\Z^2 _{odd}$, just on
$\Z^2$) will be handled in our forthcoming paper.

The infinite dimensional KAM theory with applications to Hamiltonian
PDEs has attracted great interests since  $1980$s. Starting from the
remarkable work [6,19,31], a lot of achievements have been made in
1-dimensional Hamiltonian PDEs about the existence of quasi-periodic
solutions by the methods of KAM theory. For these work, just refer
to [5,12,13,17,18,20-25,32]. But when people turn to the higher
dimensional case, the multiplicity of eigenvalues became a great
obstacle because it leads to much more complicated small divisor
conditions and measure estimates. The first breakthrough comes from
Bourgain's work [3] in $1998$. In this work, the cumbersome second
Melnikov condition is avoided due to the application of the method
of multiscale analysis, which essentially is a Nash-Moser iterative
procedure instead of Newtonian iteration being widely used in KAM
theory. Following this idea, a lot of important work has been made
in higher dimensional case (refer to [1,2,4,30]).

However, despite the advantage of avoiding the difficulty of the
second Melnikov conditions, there are also drawbacks  of multiscale
analysis methods. For example,  we couldn't see the linear stability
of the small-amplitude solutions and it couldn't show us a
description of the normal form, which is fundamental in knowing the
dynamical structure of an equation. For these reasons, KAM approach
is also expected in dealing with higher dimensional equations. The
first work comes from Geng and You [14] in 2006, which established
the KAM theorem solving higher dimensional beam equations and
nonlocal smooth Schr\"odinger equations with Fourier multiplier.
They used the ``zero-momentum condition'' to avoid the multiplicity
of eigenvalues and the regularity property to do the measure
estimate. Later in 2010, a remarkable work [8] by Eliasson and
Kuksin dealt with quite general case: higher dimensional
Schr\"odinger equations with convolutional type potential and
without ``zero-momentum condition''. To overcome the multiple
eigenvalues they studied the distribution of integer points on a
sphere and got a normal form with block-diagonal structure, and
conducted the measure estimates by developing the technique named
``Lipschitz domain''. Motivated by their method, the quasi-periodic
solutions of completely resonant Schr\"odinger equation on
2-dimension torus was developed by Geng, Xu and You [11] in 2011,
with a very elaborate construction of tangential sites. In this
paper, they defines the conception of ``T\"oplitz-Lipschitz''
condition and proved that the perturbation satisfies
``T\"oplitz-Lipschitz'' condition. Later, in [26,27] C.Procesi and
M.Procesi extended this result to higher dimensional case. For other
work about higher dimensional equation, just refer to
[7,9,10,15,16,28,29].

Let  us turn to beam equation now. In [15] Geng and You got the
quasi-periodic solutions of beam equation in high dimension with
typical constant potential and the nonlinearity is independent on
the spatial variable $x$. Recently, in [7] Eliasson, Grebert and
Kuksin got the quasi-periodic solutions of beam equation having
typical constant potential in higher dimensional case, and with an
elaborate but quite general choice of tangential sites in the sense
of probability. they allow that their normal form contain hyperbolic
terms which is cumbersome in solving homological equations.
Motivated by their work, we want to consider the completely resonant
beam equation (1.1). In our case, there are no outer parameters and
only the amplitude provides parameters. Compared with the case of
typical constant potential, although we have ``zero-momentum
condition'' here, but when doing the normal form before KAM
procedure, some terms still couldn't be eliminated because of the
loss of outer parameters.  We could only get a block-diagonal normal
form with finite dimensional block. As a consequence, our normal
form is always related to the angle variable $\theta$, so here the linear stability is not available. Compared with
[11], our convenience is that we have regularity property here and
needn't verify the complicated "T\"oeplitz-Lipschitz condition" at
each step. But except for this, our normal form structure and KAM
iteration is similar to that in [11].

Now we state the choice of tangential sites. Let $ S = \{ i_j \in
\Z^2 _{odd} : 1\leq j \leq b \} $ here $b \geq 2$. We say $S$ is
admissible if it satisfies the following conditions.

\begin{Proposition} (Structure of $S$)
$\\$
\textcircled{1}
Any three of them are not vertices of a rectangle.

\noindent \textcircled{2}
For any $n \in \Z^2 _{odd} \setminus S$, there exists at most one triple $\{ i,j,m \}$ with $i,j \in S, m \in \Z^2 _{odd} \setminus S$ such that
\begin{eqnarray*}
\left\{ \begin{array}{rcl} &&n-m+i-j=0 \\&&|n|^2 - |m|^2 + |i|^2 - |j|^2 = 0 \end{array} \right .
\end{eqnarray*}
and if it exists, we say $(n,m)$ are resonant in the first type and denote all such $n$ by ${\mathcal L}_1$.

\noindent \textcircled{3}
For any $n \in \Z^2 _{odd} \setminus S$, there exists at most one triple $\{ i,j,m \}$ with $i,j \in S, m \in \Z^2 _{odd} \setminus S$ such that
\begin{eqnarray*}
\left\{ \begin{array}{rcl} &&n+m-i-j=0 \\&&|n|^2 + |m|^2 - |i|^2 - |j|^2 = 0 \end{array} \right .
\end{eqnarray*}
and if it exists, we say $(n,m)$ are resonant in the second type and denote all such $n$ by ${\mathcal L}_2$.

\noindent \textcircled{4}
Any $n \in \Z^2 _{odd} \setminus S$ shouldn't be in ${\mathcal L}_1$ and ${\mathcal L}_2$ at the same time. It means that ${\mathcal L}_1 \cap {\mathcal L}_2 = \emptyset$.
\end{Proposition}

\noindent (Here $|\cdot|$ means $l^2$ norm)

$\\$
\indent The proof of the existence of admissible sets is postponed in
the Appendix, which is a modification of [11].

$\\$
\indent Now we could state the main theorem.

\begin{Theorem}
Let $S = (i_1, i_2, \cdots, i_b) \subseteq \Z^2 _{odd}$ be an admissible set. There exists a Cantor set $\mathcal{C}$ of positive measure, s.t. $\forall \xi = (\xi_1, \xi_2,\cdots,\xi_b) \in \mathcal{C}$, equation (1.1) admits a small-amplitude real-valued quasi-periodic solution
$$
u(t,x) = \sum_{j=1} ^b \sqrt{\xi_j} (e^{i\omega_j t}\phi_{i_j}
+e^{-i\omega_j t}\bar\phi_{i_j})+ O(|\xi|^{\frac{3}{2}})
$$
\end{Theorem}

The outline of this paper is as follows: In section 2 we state some
preliminaries and the abstract KAM theorem. In section 3 we deal
with the normal form before KAM iteration. In section 4 we conduct
one step of KAM iteration: solving homological equation and verifing
the new normal form and perturbation. In section 5 we prove uniform
convergence and get the invariant torus. In section 6 we complete
the measure estimate. The choice of tangential sites is put into the
appendix.

\section{\textbf{Preliminaries and statement of the abstract KAM theorem}}
$\\$ \indent In this section we introduce some notations and state
the abstract KAM theorem which allows the existence of some terms
dependent on $\theta$ in the normal form part.

To simplify, we only consider the subspace $\Z^2 _{odd}$ (defined in
(1.3)) instead of $\Z^2$. Given $b$ points $\{i_1,i_2,\cdots,i_b\}$
($b \geq 2$) in $\Z^2 _{odd}$, denoted by $S$, which should be an admissible set(defined in Propositon 1), and let $\Z^2 _1$ be
the complementary set of $S$ in $\Z^2 _{odd}$. Denote $z=(z_n)_{n\in
\Z^2 _1}$ with its conjugate $\bar z=(\bar z_n)_{n\in \Z^2 _1}$. We
introduce the weighted norm as follows:
\begin{eqnarray}
\| z \|_{a,\rho} = \sum_{n\in \Z^2 _1} |z_n| |n|^a e^{\rho |n|}\qquad a,\rho >0
\end{eqnarray}
\noindent Here $|n|=\sqrt{|n_1|^2 + |n_2|^2}$,$n=(n_1,n_2)\in \Z^2
_1$. Denote a neighborhood of $\T^b\times \{I=0\}\times
\{z=0\}\times\{\bar z=0\}$ by
\begin{eqnarray}
&D(r,s)=\left\{ (\theta, I, z, \bar z): |{\rm Im} \theta| < r, |I| <
s^2, \| z \|_{a,\rho} < s,\| \bar z \|_{a,\rho} < s \nonumber
\right\}
\end{eqnarray}
\noindent Here $|\cdotp|$ means the sup-norm of complex vectors.

Let $\alpha = \{\alpha_n\}_{n\in \Z^2 _1}, \beta = \{\beta_n\}_{n\in
\Z^2 _1}$, $\alpha_n, \beta_n \in \N$ with only finitely many
non-vanishing  components. Denote $z^{\alpha}\bar
z^{\beta}=\prod_{n\in \Z^2 _1}z^{\alpha_n}_n \bar z^{\beta_n}_n$ and
let \begin{eqnarray} F(\theta, I, z, \bar
z)=\sum_{k,l,\alpha,\beta}F_{kl\alpha\beta}(\xi)e^{i\left\langle
k,\theta\right\rangle}I^l z^\alpha \bar z^\beta
\end{eqnarray}
\noindent where $\xi \in {\mathcal O}\subseteq \R^b$ is the
parameter set. $k=(k_1,\cdots,k_b)\in \Z^b$ and
$l=(l_1,\cdots,l_b)\in \N^b$, $I^l = I^{l_1}_1\cdots I^{l_b}_b$.
Denote the weighted norm of $F$ by
\begin{eqnarray}
\| F \|_{D(r,s), \mathcal{O}} &=&
\sup_{\xi \in {\mathcal O}, \| z \|_{a,\rho} < s, \| \bar z \|_{a,\rho} < s}
\sum_{kl\alpha\beta}|F_{kl\alpha\beta}|_{\mathcal{O}}e^{|k|r}s^{2|l|}|z^{\alpha}||\bar z^{\beta}|   \\
|F_{kl\alpha\beta}|_{\mathcal{O}}&=&\sup_{\xi \in \mathcal{O}}\sum_{0\leq d \leq 4}|\partial_{\xi}^4 F_{kl\alpha\beta}|
\end{eqnarray}
\noindent where the derivatives with respect to $\xi$ are in the
sense of Whitney.

To a function $F$ we define its Hamiltonian vector field by
\begin{eqnarray}
X_F = (
F_I, -F_{\theta}, i\{F_{z_n}\}_{n\in \Z^2 _1}, -i\{F_{\bar z_n}\}_{n\in \Z^2 _1}
)
\end{eqnarray}
and the associated weighted norm is
\begin{eqnarray}
\|X_F\|_{D(r,s),{\mathcal O}} &:= &
\|F_I\|_{D(r,s),{\mathcal O}} + \frac{1}{s^2}\|F_{\theta}\|_{D(r,s),{\mathcal O}} \nonumber \\
 &+& \frac{1}{s}
 \left(
 \sum_{n\in \Z^2 _1} \|F_{z_n}\|_{D(r,s),{\mathcal O}}|n|^{\bar a}e^{|n|\rho} +
 \sum_{n\in \Z^2 _1} \|F_{\bar z_n}\|_{D(r,s),{\mathcal O}}|n|^{\bar a}e^{|n|\rho}
 \right)
\end{eqnarray}
\noindent where $\bar a > 0$ is a constant and we need $\bar a > a$
to measure the regularity property of the perturbation at each
iterative step.

$\\$ \indent The normal form has the following form:

\begin{eqnarray}
&H_0 = N + {\mathcal A} + {\mathcal B} + \bar {\mathcal B}  \nonumber \\
&N = \left\langle \omega(\xi),I \right\rangle + \sum\limits_{n \in \Z^2 _1} \Omega_n z_n \bar z_n  \nonumber \\
&{\mathcal A}=\sum\limits_{n\in {\mathcal L}_1} a_n (\xi) e^{i(\theta_i - \theta_j)} z_n \bar z_m  \nonumber \\
&{\mathcal B}=\sum\limits_{n\in {\mathcal L}_2} a_n (\xi) e^{-i(\theta_i + \theta_j)} z_n z_m  \nonumber \\
&\bar {\mathcal B}=\sum\limits_{n\in {\mathcal L}_2} a_n (\xi) e^{i(\theta_i + \theta_j)} \bar z_n \bar z_m  \nonumber
\end{eqnarray}
where $\xi \in \mathcal O$ is the parameter. For each $n \in
{\mathcal L}_1$ or $n \in {\mathcal L}_2$, the 3-triple $(m,i,j)$ is
uniquely determined.

For this unperturbed system, it's easy to see that it admits a
special solution $ (\theta,0,0,0)\rightarrow (\theta+\omega t,
0,0,0) $ corresponding to an invariant torus in the phase space. Our
goal is to prove that, after removing some parameters, the perturbed
system $H=H_0 + P$ still admits invariant torus provided that $
\|X_P\|_{D_{a,\rho} (r,s), {\mathcal O}} $ is sufficiently small. To
achieve this goal, we require that Hamiltonian $H$ satisfies some
conditions:

$\\$ $(\bf A1)$ Nondegeneracy: The map $\xi \rightarrow \omega
(\xi)$ is a $C^4 _W$ diffeomorphism between $\mathcal O$ and its
image ($C^4 _W$ means $C^4$ in the sense of Whitney).

\noindent $(\bf A2)$ Asymptotics of normal frequencies:
\begin{eqnarray}
\Omega_n = \varepsilon^{-p}|n|^2 + \tilde{\Omega}_n \qquad p > 0
\end{eqnarray}
here $\tilde{\Omega}_n$ is a $C^4 _W$ function of $\xi$, and
$
\tilde{\Omega}_n = O(|n|^{-\iota}) \quad \iota > 0
$

\noindent $(\bf A3)$ Melnikov conditions: Let
$$
A_n = \Omega_n \qquad n\in \Z^2 _1 \setminus ({\mathcal L}_1 \cup {\mathcal L}_2)
$$
and
$$
A_n =
\left(
\begin{array}{cc}
\Omega_n + \omega_i & a_n \\
a_m & \Omega_m + \omega_j
\end{array}
\right)
\qquad n \in {\mathcal L}_1
$$
$$
A_n =
\left(
\begin{array}{cc}
\Omega_n - \omega_i & a_n \\
\bar a_m & \Omega_m - \omega_j
\end{array}
\right)
\qquad n \in {\mathcal L}_2
$$
Then we assume that there exists $\gamma,\tau > 0$, such that
\begin{eqnarray}
&|\left\langle k,\omega \right\rangle| \geq \frac{\gamma}{|k|^{\tau}} \quad k \neq 0 \nonumber \\
&|\det(\left\langle k,\omega \right\rangle + A_n)| \geq \frac{\gamma}{|k|^{\tau}} \nonumber \\
&|\det(\left\langle k,\omega \right\rangle + A_n \otimes I_2 \pm I_2 \otimes A_{n^{'}})| \geq \frac{\gamma}{|k|^{\tau}} \quad k \neq 0 \nonumber
\end{eqnarray}

\noindent $(\bf A4)$ Boundedness:
${\mathcal A} + {\mathcal B} + \bar {\mathcal B} + P$ is real analytic in each variable $\theta,I,z,\bar z$ and Whitney smooth in $\xi$. And we have
\begin{eqnarray}
\|X_{\mathcal A}\|_{D_{a,\rho}(r,s),{\mathcal O}} +
\|X_{\mathcal B}\|_{D_{a,\rho}(r,s),{\mathcal O}} +
\|X_{\bar {\mathcal B}}\|_{D_{a,\rho}(r,s),{\mathcal O}} < 1, \qquad
\|X_P\|_{D_{a,\rho}(r,s),{\mathcal O}}<\varepsilon
\end{eqnarray}

\noindent $(\bf A5)$ Zero-momentum condition:

The normal form part ${\mathcal A} + {\mathcal B} + \bar {\mathcal
B} + P$ satisfy the following condition:
$$
{\mathcal A} + {\mathcal B} + \bar {\mathcal B} + P= \sum\limits_{k
\in \Z^b, l \in \N^b, \alpha, \beta} ({\mathcal A} + {\mathcal B} +
\bar {\mathcal B} + P)_{kl\alpha \beta} (\xi) e^{i\left\langle
k,\theta \right\rangle}I^l z^{\alpha} {\bar z}^{\beta}
$$
we have
$$
({\mathcal A} + {\mathcal B} + \bar {\mathcal B} + P)_{kl\alpha
\beta} \neq 0 \Longrightarrow \sum_{j=1} ^b k_j i_j +
\sum\limits_{n\in \Z^2 _1}(\alpha_n - \beta_n)n = 0
$$

Now we state our abstract KAM theorem, and as a corollary, we get Theorem $\bf 1$.

\begin{Theorem}
Assume that the Hamiltonian $H=N+ {\mathcal A} + {\mathcal B} + \bar
{\mathcal B} + P$ satisfies condition $\bf (A1)-(A5)$. Let $\gamma >
0$ be sufficiently small, then there exists $\varepsilon > 0$ and
$a,\rho > 0$ such that if $\|X_P\|_{D_{a,\rho}(r,s),\mathcal O} <
\varepsilon$, the following holds: There exists a Cantor subset
${\mathcal O}_{\gamma} \subseteq {\mathcal O}$ with $meas({\mathcal
O}\setminus {\mathcal O}_{\gamma}) = O(\gamma^{\varsigma})$
($\varsigma$ is a positive constant) and two maps which are analytic
in $\theta$ and $C_W ^4$ in $\xi$.
$$
\Phi:\T^b \times {\mathcal O}_{\gamma} \rightarrow D_{a,\rho}(r,s),\qquad \tilde{\omega}: {\mathcal O}_{\gamma}\rightarrow \R^b
$$
where $\Phi$ is $\frac{\varepsilon}{\gamma^{16}}$-close to the
trivial embedding $\Phi_0 : \T^b \times {\mathcal O} \rightarrow
\T^b \times \{0,0,0\}$ and $\tilde{\omega}$ is $\varepsilon$-close
to the unperturbed frequency $\omega$, such that $\forall \xi \in
{\mathcal O}_{\gamma}$ and $\theta \in \T^b$, the curve
$t\rightarrow \Psi(\theta+ \tilde{\omega}t, \xi)$ is a
quasi-periodic solution of the Hamiltonian equation governed by
$H=N+{\mathcal A} + {\mathcal B} + \bar {\mathcal B} + P$.

\end{Theorem}

\section{\textbf{Normal Form}}
$\\$
\indent Consider the equation (1.1). The linear operator $-\Delta$ has eigenvalues $\lambda_n = |n|^2$and corresponding eigenfunctions $\phi_n = \frac{1}{2\pi}e^{i\left\langle n,x \right\rangle}$. By scaling $u\rightarrow \varepsilon^{\frac{1}{2}}u$, (1.1)becomes
\begin{eqnarray}
u_{tt}+\Delta^2 u +\varepsilon u^3 =0
\end{eqnarray}
Now introduce $v = u_t$ and (3.1) is turned into
\begin{eqnarray}
&&u_t = v \nonumber\\
&&v_t = -\Delta^2 u -\varepsilon u^3
\end{eqnarray}
Let $q=\frac{1}{\sqrt 2}((-\Delta)^{\frac{1}{2}}u - i(-\Delta)^{-\frac{1}{2}}v)$ and (3.2) becomes
\begin{eqnarray}
-iq_t=-\Delta q + \varepsilon\frac{1}{\sqrt 2}(-\Delta)^{-\frac{1}{2}}\left((-\Delta)^{-\frac{1}{2}}(\frac{q+\bar q}{\sqrt 2})\right)^3
\end{eqnarray}

Write it in the form of Hamiltonian equation $q_t = i\frac{\partial H}{\partial \bar q}$ and we get the Hamiltonian
\begin{eqnarray}
H=\frac{1}{2}\left\langle -\Delta q, q \right\rangle+\frac{1}{4}\varepsilon\int_{\T^2}\left( (-\Delta)^{-\frac{1}{2}}(q+\bar q) \right)^4 dx
\end{eqnarray}
where $\left\langle \cdot,\cdot \right\rangle$ is the inner product
in $L^2 (\T^2)$. Notice that in $\Z^2 _{odd}$ the origin is
avoided so $(-\Delta)^{-\frac{1}{2}}$ is well defined. (That is why we
use it instead of the whole $\Z^2$) Now expand $q$ into Fourier
series
\begin{eqnarray}
q=\sum_{n\in \Z^2 _{odd}}q_n \phi_n
\end{eqnarray}
so the Hamiltonian becomes (justify $\varepsilon$ if necessary)
\begin{eqnarray}
H=\sum_{n\in \Z^2 _{odd}}\lambda_n |q_n|^2 &+& \varepsilon \sum_{i+j+k+l=0}\frac{1}{\sqrt{\lambda_i \lambda_j \lambda_k \lambda_l}}
\left( q_i q_j q_k q_l +\bar{q_i}\bar{q_j}\bar{q_k}\bar{q_l} \right) \nonumber \\
&+&4\varepsilon \sum_{i+j+k-l=0}\frac{1}{\sqrt{\lambda_i \lambda_j \lambda_k \lambda_l}}
\left(q_i q_j q_k \bar q_l + \bar q_i \bar q_j \bar q_k q_l \right) \nonumber \\
&+&6\epsilon \sum_{i+j-k-l=0}\frac{1}{\sqrt{\lambda_i \lambda_j \lambda_k \lambda_l}}
\left(q_i q_j \bar{q_k}\bar{q_l}\right)
\end{eqnarray}
Now we state the normal form theorem of $H$.
\newtheorem{proposition}{Propsition}[section]
\begin{proposition}
Let $S$ be admissible. For Hamiltonian function (3.6), there exists a symplectic transformation $\Phi$ satisfying
\begin{eqnarray}
H\circ \Phi = \left\langle \omega, I \right\rangle + \sum_{n\in \Z^2 _1}\Omega_n z_n \bar z_n
+{\mathcal A} + {\mathcal B} + \bar{\mathcal B} + P
\end{eqnarray}
where
\begin{eqnarray}
\omega_i &=& \varepsilon^{-4}\lambda_i + \frac{2}{\lambda_i ^2}\xi_i + 4\sum_{j\in S,j\neq i}\frac{1}{\lambda_i \lambda_j}\xi_j \quad i\in S  \\
\Omega_n &=& \varepsilon^{-4}\lambda_n + 4\sum_{j\in S}\frac{1}{\lambda_j \lambda_n}\xi_j \quad n\in \Z^2 _1
\end{eqnarray}
and
\begin{eqnarray}
\mathcal A &=&4\sum_{n\in {\mathcal L}_1}\frac{\sqrt{\xi_i \xi_j}}{\sqrt{\lambda_i \lambda_j \lambda_n \lambda_m}}e^{i(\theta_i - \theta_j)}z_n \bar z_m \\
\mathcal B &=&4\sum_{n\in {\mathcal L}_2}\frac{\sqrt{\xi_i \xi_j}}{\sqrt{\lambda_i \lambda_j \lambda_n \lambda_m}}e^{i(-\theta_i - \theta_j)}z_n z_m \\
\bar {\mathcal B} &=& 4\sum_{n\in {\mathcal L}_2}\frac{\sqrt{\xi_i \xi_j}}{\sqrt{\lambda_i \lambda_j \lambda_n \lambda_m}}e^{i(\theta_i + \theta_j)}\bar z_n \bar z_m
\end{eqnarray}
\begin{eqnarray}
P&=& O (\varepsilon^2 |I|^2 + \varepsilon^2 |I|\| z \| _{a,\rho}^2 +
\varepsilon |\xi|^{\frac{1}{2}}\| z \| _{a,\rho}^3
 + \varepsilon^2 \| z \| _{a,\rho}^4 + \varepsilon^2 |\xi|^3 \nonumber \\
&+& \varepsilon^3 |\xi|^{\frac{5}{2}}\| z \|_{a,\rho}+\varepsilon^4
|\xi|^2 \| z \|_{a,\rho} ^2 +\varepsilon^5 |\xi|^{\frac{3}{2}}\| z
\|_{a,\rho}^3 )
\end{eqnarray}
\end{proposition}
\noindent $\bf Proof:$
We construct a Hamiltonian function $F$ to induce $\Phi = X^1 _F$ which is the time-1 map of $F$. For convenience, we define three sets as below:
\begin{eqnarray}
S_1 = \{(i,j,n,m): &&\textcircled{1}: i-j+n-m=0 \nonumber \\
&&\textcircled{2}: |i|^2 - |j|^2 + |n|^2 - |m|^2 \neq 0 \nonumber \\
&&\textcircled{3}: \#\{i,j,n,m\}\cap S \geq 2  \}
\end{eqnarray}
and similarly
\begin{eqnarray}
S_2 = \{(i,j,n,m): &&\textcircled{1}: i+j+n+m=0 \nonumber \\
&&\textcircled{2}: |i|^2 + |j|^2 + |n|^2 + |m|^2 \neq 0 \nonumber \\
&&\textcircled{3}: \#\{i,j,n,m\}\cap S \geq 2  \}
\end{eqnarray}
\begin{eqnarray}
S_3 =\{(i,j,n,m): &&\textcircled{1}: i+j+n-m=0 \nonumber \\
&&\textcircled{2}: |i|^2 + |j|^2 + |n|^2 - |m|^2 \neq 0 \nonumber \\
&&\textcircled{3}: \#\{i,j,n,m\}\cap S \geq 2  \}
\end{eqnarray}

we define $F$ as
\begin{eqnarray}
F &=& \sum_{S_1}\frac{i\varepsilon}{\lambda_i - \lambda_j + \lambda_n - \lambda_m}q_i\bar q_j q_n \bar q_m \nonumber \\
  &+& \sum_{S_2}\frac{i\varepsilon}{6(\lambda_i + \lambda_j + \lambda_n + \lambda_m)}(q_iq_jq_nq_m - \bar q_i \bar q_j \bar q_n \bar q_m) \nonumber \\
  &+& \sum_{S_3}\frac{2i\varepsilon}{3(\lambda_i + \lambda_j + \lambda_n - \lambda_m)}(q_iq_jq_n\bar q_m - \bar q_i \bar q_j \bar q_n q_m)
\end{eqnarray}
(3.6) is put into (set $z_n = q_n,\bar z_n = \bar q_n, n \notin S$)
\begin{eqnarray}
H\circ \Phi &=& \sum_{n\in S}\lambda_n |q_n|^2 + \sum_{n\notin S}\lambda_n |z_n|^2 + \varepsilon \sum_{n\in S}\frac{1}{\lambda_n ^2}|q_n|^4 \nonumber \\
&+& 4\varepsilon \sum_{i,j\in S, i\neq j}\frac{1}{\lambda_i \lambda_j}|q_i|^2|q_j|^2 + 4\varepsilon \sum_{i\in S, n\notin S}\frac{1}{\lambda_i \lambda_n}|q_i|^2|z_n|^2 \nonumber \\
&+& 4\varepsilon \sum_{n\in {\mathcal L}_1}\frac{1}{\sqrt{\lambda_i \lambda_j \lambda_n \lambda_m}}q_i\bar q_j z_n{\bar z}_m + 4\varepsilon \sum_{n\in {\mathcal L}_2}\frac{1}{\sqrt{\lambda_i \lambda_j \lambda_n \lambda_m}}(q_i q_j\bar z_n\bar z_m + \bar q_i \bar q_j z_n z_m) \nonumber \\
&+& O\left(\varepsilon |q|||z||^3 _{a,\rho} +\varepsilon ||z||^4
_{a,\rho}+\varepsilon^2 |q|^6 + \varepsilon^2 |q|^5 \|z\|^3
_{a,\rho} + \varepsilon^2 |q|^4 ||z||^2 _{a,\rho} + \varepsilon^2
|q|^3 ||z||^3 _{a,\rho}  \right) \nonumber
\end{eqnarray}

Here we need to state a fact: For four points $n,m,i,j \in \Z^2 _{odd}$ , it could never satisfy
$|n|^2 + |m|^2 + |i|^2 - |j|^2 = 0$. If not, we assume $n=(n_1, n_2),m=(m_1,m_2),i=(i_1,i_2),j=(j_1,j_2)$ and in each one the first component is odd and the second component is even. Then we have
$$
|n_1|^2 + |m_1|^2 + |i_1|^2 - |j_1|^2 = -
(
|n_2|^2 + |m_2|^2 + |i_2|^2 - |j_2|^2
)
$$
The right one can be divided by 4 but the left one couldn't, which is a contradiction. By this fact we know that the set
\begin{eqnarray}
\{(i,j,n,m) \in ( \Z^2 _{odd} )^4: &&\textcircled{1}: i+j+n-m=0 \nonumber \\
&&\textcircled{2}: |i|^2 + |j|^2 + |n|^2 - |m|^2 = 0 \}
\end{eqnarray}
is empty.

Introduce the action-angle variable in the tangential sites:
\begin{eqnarray}
q_j=\sqrt{I_j + \xi_j}e^{i\theta_j},\quad \bar q_j=\sqrt{I_j + \xi_j}e^{-i\theta_j}
\end{eqnarray}
so we have

\begin{eqnarray}
H \circ \Phi &=& \sum_{i \in S}\lambda_i (I_i + \xi_i) + \sum_{n\notin S}\lambda_n |z_n|^2 + \varepsilon
\sum_{i\in S}\frac{1}{\lambda_i ^2}(I_i + \xi_i)^2 \nonumber \\
&+& 4\varepsilon \sum_{i,j\in S, i\neq j}\frac{1}{\lambda_i \lambda_j}(I_i + \xi_i)(I_j + \xi_j) +
    4\varepsilon \sum_{i\in S, n\notin S}\frac{1}{\lambda_i \lambda_n}(I_i + \xi_i)|z_n|^2 \nonumber \\
&+& 4\varepsilon \sum_{n\in {\mathcal L}_1}\frac{1}{\sqrt{\lambda_i \lambda_j \lambda_n \lambda_m}}
    \sqrt{(I_i + \xi_i)(I_j + \xi_j)}e^{i(\theta_i - \theta_j)}z_n \bar z_m \nonumber \\
&+& 4\varepsilon \sum_{n\in {\mathcal L}_2}\frac{1}{\sqrt{\lambda_i \lambda_j \lambda_n \lambda_m}}
    \sqrt{(I_i + \xi_i)(I_j + \xi_j)}e^{i(-\theta_i - \theta_j)}z_n z_m \nonumber \\
&+& 4\varepsilon \sum_{n\in {\mathcal L}_2}\frac{1}{\sqrt{\lambda_i \lambda_j \lambda_n \lambda_m}}
    \sqrt{(I_i + \xi_i)(I_j + \xi_j)}e^{i(\theta_i +\theta_j)}\bar z_n \bar z_m \nonumber \\
&+& O(\varepsilon |\xi|^{\frac{1}{2}}\|z\|_{a,\rho}^3 + \varepsilon \|z\|_{a,\rho}^4 + \varepsilon^2 |\xi|^3 +
\varepsilon^2 |\xi|^{\frac{5}{2}}\|z\|_{a,\rho} + \varepsilon^2 |\xi|^2 \|z\|_{a,\rho}^2
+ \varepsilon^2 |\xi|^{\frac{3}{2}}\|z\|_{a,\rho}^3)\nonumber \\
&=& \sum_{i\in S}(\lambda_i + \frac{2\varepsilon}{\lambda_i ^2}\xi_i + 4\varepsilon \sum_{j\in S, i\neq j}\frac{1}{\lambda_i \lambda_j}\xi_j)\omega_i \nonumber \\
&+& \sum_{n\notin S}(\lambda_n + 4\varepsilon \sum_{i\in S}\frac{1}{\lambda_i \lambda_n}\xi_i)|z_n|^2 \nonumber \\
&+& 4\varepsilon \sum_{n\in {\mathcal L}_1}\frac{1}{\sqrt{\lambda_i \lambda_j \lambda_n \lambda_m}}
    \sqrt{\xi_i \xi_j}e^{i(\theta_i - \theta_j)}z_n \bar z_m \nonumber \\
&+& 4\varepsilon \sum_{n\in {\mathcal L}_2}\frac{1}{\sqrt{\lambda_i \lambda_j \lambda_n \lambda_m}}
    \sqrt{\xi_i \xi_j}e^{i(-\theta_i - \theta_j)}z_n z_m \nonumber \\
&+& 4\varepsilon \sum_{n\in {\mathcal L}_2}\frac{1}{\sqrt{\lambda_i \lambda_j \lambda_n \lambda_m}}
    \sqrt{\xi_i \xi_j}e^{i(\theta_i + \theta_j)}\bar z_n \bar z_m \nonumber \\
&+& O(\varepsilon |I|^2 + \varepsilon |I| \|z\|_{a,\rho}^2 +
\varepsilon |\xi|^{\frac{1}{2}}\|z\|_{a,\rho}^3 + \varepsilon \|z\|_{a,\rho}^4 +
\varepsilon^2 |\xi|^3 \nonumber \\
&+& \varepsilon^2 |\xi|^{\frac{5}{2}}\|z\|_{a,\rho} + \varepsilon^2 |\xi|^2
\|z\|_{a,\rho}^2 + \varepsilon^2 |\xi|^{\frac{3}{2}}\|z\|_{a,\rho}^3)
\nonumber
\end{eqnarray}
By scaling in variables:
$$
\xi \rightarrow \varepsilon^3 \xi, \quad I\rightarrow \varepsilon^5 I,\quad z\rightarrow \varepsilon^{\frac{5}{2}}z,
\quad \bar z\rightarrow \varepsilon^{\frac{5}{2}}\bar z
$$
and scale time $t\rightarrow \varepsilon^9 t$ we get the Hamiltonian
function as follows:
\begin{eqnarray}
H=\left\langle \omega, I\right\rangle + \left\langle \Omega z, z\right\rangle + {\mathcal A} + {\mathcal B} + \bar{\mathcal B} + P
\end{eqnarray}
where
\begin{eqnarray}
\omega_i &=& \varepsilon^{-4}\lambda_i + \frac{2}{\lambda_i ^2}\xi_i + 4\sum_{j\in S,j\neq i}\frac{1}{\lambda_i \lambda_j}\xi_j  \\
\Omega_n &=& \varepsilon^{-4}\lambda_n + 4\sum_{j\in S}\frac{1}{\lambda_j \lambda_n}\xi_j  \\
{\mathcal A} &=& 4\sum_{n\in {\mathcal L}_1}\frac{\sqrt{\xi_i \xi_j}}{\sqrt{\lambda_i \lambda_j \lambda_n \lambda_m}}
               e^{i(\theta_i - \theta_j)}z_n \bar z_m  \\
{\mathcal B} &=& 4\sum_{n\in {\mathcal L}_2}\frac{\sqrt{\xi_i \xi_j}}{\sqrt{\lambda_i \lambda_j \lambda_n \lambda_m}}
               e^{i(-\theta_i - \theta_j)}z_n z_m  \\
\bar{\mathcal B} &=& 4\sum_{n\in {\mathcal L}_2}\frac{\sqrt{\xi_i \xi_j}}{\sqrt{\lambda_i \lambda_j \lambda_n \lambda_m}}
               e^{i(\theta_i + \theta_j)}\bar z_n \bar z_m  \\
P &=& O( \varepsilon^2 |I|^2 + \varepsilon^2 |I| \|z\|_{a,\rho}^2 + \varepsilon |\xi|^{\frac{1}{2}}\|z\|_{a,\rho}^3 +
\varepsilon^2\|z\|_{a,\rho}^4 +\varepsilon^2 |\xi|^3 + \varepsilon^3 |\xi|^{\frac{5}{2}}\|z\|_{a,\rho} \nonumber \\
&+& \varepsilon^4 |\xi|^2 \|z\|_{a,\rho}^2 + \varepsilon^5
|\xi|^{\frac{3}{2}}\|z\|_{a,\rho} ^3 )
\end{eqnarray}

Now we verify that the normal form $(3.7)-(3.13)$ satisfy condition $\bf (A1)-(A5)$.

$Verifying \bf (A1)$: By (3.8) we get
\begin{eqnarray}
\frac{\partial \omega}{\partial \xi}= (a_{ij})_{1\leq i,j \leq b}
\end{eqnarray}
where $a_{ij}=\frac{2}{\lambda_i ^2}$ if $i=j$ and
$a_{ij}=\frac{4}{\lambda_i \lambda_j}$ if $i \neq j$. It's easy to
see that this matrix is non--degenerate.

$Verifying \bf (A2)$: By $(3.9)$, just take $p=4,\iota=2$.

$Verifying \bf (A3)$: Recall the definition in condition $(A3)$, we only verify the most complicated case:
\begin{eqnarray}
\det (\left\langle k,\omega \right\rangle + A_n \otimes I_2 - I_2 \otimes A_{n^{'}})
\end{eqnarray}
where $n,n^{'} \in {\mathcal L}_1 \cup {\mathcal L}_2$. We verify
two facts: (3.28) is a polynomial of parameter $\xi$ with degree 4
and it couldn't be equivalently zero. For the former one, notice
that $ \lambda I + A\otimes I_2 - I_2 \otimes B = (\lambda I +
A)\otimes I-I \otimes B $ (here $|\cdot|$ means determinant) and
using the formula
$$
|A\otimes I \pm I\otimes B|=
(|A|-|B|)^2 +|A|(tr(B))^2 + |B|(tr(A))^2 \pm (|A|+|B|)tr(A)tr(B)
$$
then we get it. For the latter one, it's the same as that in [11].
By this, we could get
$$
|\partial_{\xi} ^4
\left(
\det (
\left\langle k,\omega \right\rangle +A_n \otimes I_2 \pm I_2 \otimes A_{n^{'}}
)
\right)| > c |k|
$$
So by excluding parameters with measure $O(\gamma^{\frac{1}{4}})$, we have
$$
|\left(
\det (
\left\langle k,\omega \right\rangle +A_n \otimes I_2 \pm I_2 \otimes A_{n^{'}}
)
\right)| > \frac{\gamma}{|k|^{\tau}} \quad k \neq 0
$$

For the verification of $\bf (A4)$ and $\bf (A5)$, just refer to [14].

$\\$
\section{\textbf{KAM Iteration}}
$\\$
 \indent We prove Theorem $\bf 2$ by a KAM iteration which involves an infinite sequence of change of variables. Each step of KAM iteration makes the perturbation smaller than the previous step at the cost of excluding a small set of parameters. We have to prove the convergence of the iteration and estimate the measure of the excluded set after infinite KAM steps.

At the $\nu$-step of the KAM iteration, we consider a Hamiltonian vector field with
$$
H_{\nu} = N_{\nu} + {\mathcal A}_{\nu} + {\mathcal B}_{\nu} + \bar {\mathcal B}_{\nu} + P_{\nu} =
\left\langle\omega_{\nu},I \right\rangle + \sum_{n\in \Z^{2}_1}\Omega^{\nu}_n z_n \bar z_n + {\mathcal A}_{\nu} + {\mathcal B}_{\nu} + \bar {\mathcal B}_{\nu} + P_{\nu}
$$
where $A_{\nu} + B_{\nu} + \bar B_{\nu} + P_{\nu}$ is defined in $D(r_{\nu}, s_{\nu})\times {\mathcal O}_{\nu -1}$

We construct a map
$$
\Phi_{\nu}:D(r_{\nu+1},s_{\nu+1})\times {\mathcal O}_{\nu} \rightarrow D(r_{\nu},s_{\nu})\times {\mathcal O}_{\nu-1}
$$
so that the vector field $X_{H_{\nu}\circ \Phi_{\nu}}$ defined on $D(r_{\nu+1},s_{\nu+1})$satisfies
$$
\| X_{P_{\nu+1}} \|_{ D(r_{\nu+1},s_{\nu+1}), {\mathcal O}_{\nu} }=\| X_{H_{\nu}\circ \Phi_{\nu}} -
X_{N_{\nu+1}+{\mathcal A}_{\nu+1}+{\mathcal B}_{\nu+1}+\bar{\mathcal B}_{\nu+1}} \|_{ D(r_{\nu+1},s_{\nu+1}), {\mathcal O}_{\nu} } \leq \varepsilon_{\nu} ^{\kappa}, \quad \kappa > 1
$$
and the new Hamiltonian still satisfies $\bf (A1)-(A5)$.

\indent To simplify notations, in the following text, the quantities without subscripts refer to quantities at the $\nu$ th step, while the quantities with subscripts $+$ denote the corresponding quantities at the $(\nu+1)$ th step. Let's consider the Hamiltonian defined in $D(r,s)\times \mathcal O$:
\begin{eqnarray}
H&=& N+{\mathcal A} + {\mathcal B} + \bar{\mathcal B} + P \nonumber \\
 &=& e+\left\langle \omega(\xi), I \right\rangle + \sum_{n\in \Z^{2}_1}\Omega_n z_n \bar z_n +
     {\mathcal A} + {\mathcal B} + \bar {\mathcal B} +
     P(\theta, I, z, \bar z, \xi, \varepsilon)
\end{eqnarray}
We assume that for $\xi \in \mathcal O$, one has
\begin{eqnarray}
&&|\left\langle k,\omega \right\rangle| \geq \frac{\gamma}{|k|^{\tau}} \quad k \neq 0 \nonumber \\
&&|\det( \left\langle k,\omega \right\rangle I + A_n )| \geq  \frac{\gamma}{|k|^{\tau}} \nonumber \\
&&|\det( \left\langle k,\omega \right\rangle I \pm A_n \otimes I_2
\pm I_2 \otimes A_{n^{\prime}} )| \geq  \frac{\gamma}{|k|^{\tau}}
\quad k\neq 0
\end{eqnarray}
where $A_n = \Omega_n$ for $n\in {\Z^2 _1} \setminus ({\mathcal L}_1 \cup {\mathcal L}_2)$
and
$$
A_n =
\left(
\begin{array}{cc}
\Omega_n + \omega_i & a_n \\
a_m & \Omega_m + \omega_j
\end{array}
\right)
\qquad n \in {\mathcal L}_1
$$
$$
A_n =
\left(
\begin{array}{cc}
\Omega_n - \omega_i & a_n \\
\bar a_m & \Omega_m - \omega_j
\end{array}
\right)
\qquad n \in {\mathcal L}_2
$$
where $(n,m)$ are resonant pairs and $(i,j)$ is uniquely determined by $(n,m)$.

Expand $P$ into Fourier-Taylor series
$
P=\sum_{k,l,\alpha,\beta}P_{kl\alpha\beta}I^l e^{i\left\langle k,\theta \right\rangle}z^{\alpha}\bar z^{\beta}
$
and by $(A5)$ we get that
\begin{eqnarray}
P_{kl\alpha\beta}=0 \quad if \quad
\sum_{1\leq j \leq b}k_j i_j + \sum_{n\in \Z^2 _1}(\alpha_n - \beta_n)n \neq 0
\end{eqnarray}

We now let $0<r_+<r$ and define
\begin{eqnarray}
s_+ = \frac{1}{4}s\varepsilon^{\frac{1}{3}}, \quad \varepsilon_+ = c\gamma^{-16}(r-r_+)^{-c}\varepsilon^{\frac{4}{3}}
\end{eqnarray}
Here and later, the letter $c$ denotes suitable(possible different) constants independent of the iteration steps.

Now we describe how to construct a set $\mathcal O_+ \subseteq \mathcal O$ and a change of variables
$\Phi : D_+ \times {\mathcal O}_+ = D(r_+, s_+)\times {\mathcal O}_+ \rightarrow D(r,s)\times {\mathcal O}$ such that the transformed Hamiltonian $H_+ = H\circ \Phi = N + {\mathcal A}_+ + {\mathcal B}_+ + \bar{\mathcal B}_+ + P_+$ satisfies assumptions $\bf (A1)-(A5)$ with new parameters $\varepsilon_+, r_+, s_+$ and with $\xi \in \mathcal O_+$.

\subsection{\textbf{Homological Equation}}

$\\$
\indent Expand $P$ into Fourier-Taylor series
\begin{eqnarray}
P=\sum_{k,l,\alpha,\beta}P_{kl\alpha\beta}I^l e^{i\left\langle k,\theta \right\rangle}z^{\alpha}\bar z^{\beta}
\end{eqnarray}
where $k\in \Z^b, l\in \N^b$and the multi-indices $\alpha$,$\beta$
run over the set of all infinite dimensional vectors $\alpha =
(\cdots ,\alpha_n,\cdots)_{n\in \Z^2 _1}$,  $\beta = (\cdots
,\beta_n,\cdots)_{n\in \Z^2 _1}$ with finitely many nonzero
components of positive integers. And by $(A5)$ we get that
\begin{eqnarray}
P_{kl\alpha\beta}=0 \quad if \quad
\sum_{1\leq j \leq b}k_j i_j + \sum_{n\in \Z^2 _1}(\alpha_n - \beta_n)n \neq 0
\end{eqnarray}

Consider its quadratic truncation $R$:
\begin{eqnarray}
R(\theta, I, z, \bar z)
&=& R_0 + R_1 + R_2 \nonumber \\
&=& \sum_{k,|l|\leq 1}P_{kl00}e^{i\left\langle k,\theta \right\rangle}I^l \nonumber \\
&+& \sum_{k,n}( P^{k10}_n z_n + P^{k01}_n \bar z_n  )e^{i\left\langle k,\theta \right\rangle} \nonumber \\
&+& \sum_{k,n,m}(P^{k20}_{nm}z_n z_m + P^{k11}_{nm}z_n \bar z_m + P^{k02}_{nm}\bar z_n \bar z_m)e^{i\left\langle k,\theta \right\rangle}
\end{eqnarray}
where $P^{k10}_n = P_{kl\alpha \beta}$ with $\alpha = e_n$,$\beta = 0$, $P^{k01}_n = P_{kl\alpha \beta}$ with $\alpha = 0$,$\beta = e_n$, here $e_n$ denotes the vector with the $n$ th component being $1$ and the other components being zero. Similarly, $P^{k20}_{nm}=P_{kl\alpha \beta}$ with $\alpha = e_n + e_m$, $\beta = 0$; $P^{k11}_{nm}=P_{kl\alpha \beta}$ with $\alpha = e_n$,$\beta = e_m$; $P^{k02}_{nm}=P_{kl\alpha \beta}$ with $\alpha =0$, $\beta = e_n + e_m$.

Rewrite $H$ as $H= N+ {\mathcal A} + {\mathcal B} + \bar{\mathcal B}
+ R + (P - R)$. Due to the choice of $s_+ \ll s$ and the definition
of the norm, it follows  immediately
\begin{eqnarray}
\| X_R \|_{D(r,s),\mathcal O}\leq \| X_P \|_{D(r,s),\mathcal O}\leq \varepsilon
\end{eqnarray}
and in $D(r,s_+)$
\begin{eqnarray}
\| X_{P-R} \|_{D(r,s_+)}\leq \varepsilon_+
\end{eqnarray}

In the following, we will construct a Hamiltonian function $F$ satisfying $(A5)$ and with the same form of $R$ defined in $D_+ = D(r_+ , s_+)$ such that the time one map $X^1 _F$ of the Hamiltonian vector field $X_F$ defines a map from $D_+$ to $D$ and puts $H$ into $H_+$. Precisely, one has
\begin{eqnarray}
H\circ X^1 _F
&=& (N+{\mathcal A} + {\mathcal B} + \bar{\mathcal B}+ R)\circ X^1 _F + (P-R)\circ X^1 _F \nonumber \\
&=& N+{\mathcal A} + {\mathcal B} + \bar{\mathcal B}  \nonumber \\
&+& \{N+ {\mathcal A} + {\mathcal B} + \bar{\mathcal B}, F \} + R  \\
&+& \int_0 ^1 (1-t)\{ \{N+{\mathcal A} + {\mathcal B} + \bar{\mathcal B}, F \}, F \}\circ X^t _F dt \\
&+& \int_0 ^1 \{R, F\}\circ X^t _F dt + (P-R)\circ X^1 _F
\end{eqnarray}
So we get the linearized homological equation:
\begin{eqnarray}
\{N+ {\mathcal A} + {\mathcal B} + \bar{\mathcal B}, F \}+ R = {\hat N} + \hat{\mathcal A}+ \hat{\mathcal B} + \hat{\bar{\mathcal B}}
\end{eqnarray}
where
\begin{eqnarray}
\hat N &=& P_{0000} + \left\langle \hat{\omega}, I \right\rangle + \sum_{n\in \Z^2 _1}P^{011}_{nn} z_n \bar z_n
\quad \hat{\omega}=(P_{0100})_{|l|=1}\\
\hat{\mathcal A} &=& \sum_{n\in {\mathcal L}_1}P^{11}_{nm}e^{i(\theta_i-\theta_j)}z_n \bar z_m \\
\hat{\mathcal B} &=& \sum_{n\in {\mathcal L}_2}P^{20}_{nm}e^{-i(\theta_i+\theta_j)}z_n z_m \\
\hat{\bar{\mathcal B}} &=& \sum_{n\in {\mathcal L}_2}P^{02}_{nm}e^{i(\theta_i+\theta_j)}\bar z_n\bar z_m
\end{eqnarray}
We define
$
N_+ = N+\hat N, {\mathcal A}_+ = {\mathcal A} + \hat{\mathcal A}, {\mathcal B}_+ = {\mathcal B} + \hat{\mathcal B},
\bar{\mathcal B_+} = \bar{\mathcal B} + \hat{\bar{\mathcal B}}
$
and
\begin{eqnarray}
P_+ &=& \int_0 ^1 (1-t)\{ \{N+{\mathcal A} + {\mathcal B} + \bar{\mathcal B}, F \}, F \}\circ X^t _F dt  \nonumber\\
&+& \int_0 ^1 \{R, F\}\circ X^t _F dt + (P-R)\circ X^1 _F
\end{eqnarray}
We construct the Hamiltonian function $F$ as below:
\begin{eqnarray}
F(\theta, I ,z , \bar z)
&=& F^0 + F^1 + F^{10} + F^{01} + F^{20} + F^{11} + F^{02} \nonumber \\
&=& F^0 (\theta) +\left\langle F^1(\theta), I\right\rangle + \left\langle F^{10}(\theta), z\right\rangle
    +\left\langle F^{01}(\theta), \bar z\right\rangle \nonumber \\
&+& \left\langle F^{20}(\theta)z, z\right\rangle + \left\langle
F^{11}(\theta)z, \bar z\right\rangle
    + \left\langle F^{02}(\theta)\bar z, \bar z\right\rangle
\end{eqnarray}
Now $(4.13)$ turns to
\begin{eqnarray}
\{N, F^0 + F^1\}+ R^0 + R^1 - P_{0000} - \left\langle \hat{\omega}, I\right\rangle &=& 0  \\
\{N+ {\mathcal A} + {\mathcal B} + \bar{\mathcal B}, F^{10}+ F^{01}\}+ R^{10} + R^{01} &=& 0
\end{eqnarray}
and the most complicated
\begin{eqnarray}
&&\{N +{\mathcal A} + {\mathcal B} + \bar{\mathcal B}, F^{20}+F^{11}+F^{02}\}+R^{20}+R^{11}+R^{02} \nonumber \\
&& = \sum_{n\in \Z^2 _1}P^{011}_{nn}z_n \bar z_n + \hat{\mathcal A} + \hat{\mathcal B} + \hat{\bar{\mathcal B}}
\end{eqnarray}
$\\$
$\bf Solving (4.20)$:\quad
$
F^j (\theta)=\sum\limits_{k\neq 0}F^j _k e^{i\left\langle k,\theta\right\rangle}, \quad j=0,1.
$
By comparing the Fourier coefficients we get
$$
F^j _k = \frac{-i}{\left\langle k,\omega\right\rangle}P^j _k,\quad j=0,1, k\neq 0
$$
and according to assumption $(4.2)$ we get
$$
|\left\langle k,\omega(\xi)\right\rangle|\geq
\frac{\gamma}{|k|^{\tau}}, k\neq 0, \xi \in \mathcal O
$$
so one has the estimate
\begin{eqnarray}
|F^j _k|_{\mathcal O} \leq \gamma^{-16}k^{16\tau + 16}|P^j _k|_{\mathcal O}
\end{eqnarray}

\noindent $\bf Solving (4.21)$:
We decompose this part into three cases:

\noindent $(1)$ If $n\in \Z^2 _1 \setminus \{{\mathcal L}_1 \cup {\mathcal L}_2\}$, one has
\begin{eqnarray}
(\left\langle k,\omega\right\rangle + \Omega_n)F^{10}_{k,n} = -iR^{10}_{k,n} \nonumber \\
(\left\langle k,\omega\right\rangle - \Omega_n)F^{10}_{k,n} = -iR^{01}_{k,n}
\end{eqnarray}
\noindent $(2)$ If $n\in {\mathcal L}_1$ and the corresponding resonant group $£¨m,i,j£©$, one has
\begin{eqnarray}
(\left\langle k+e_i,\omega\right\rangle + \Omega_n)F^{10}_{k+e_i,n} + a_n F^{10}_{k+e_j,m} = -iR^{10}_{k+e_i,n} \nonumber\\
(\left\langle k+e_j,\omega\right\rangle + \Omega_m)F^{10}_{k+e_j,m} + a_m F^{10}_{k+e_i,n} = -iR^{10}_{k+e_j,m}
\end{eqnarray}
\noindent $(3)$ If $n\in {\mathcal L}_2$, one has
\begin{eqnarray}
(\left\langle k-e_i,\omega\right\rangle + \Omega_n)F^{10}_{k-e_i,n} - a_n F^{01}_{k+e_j,m} =-iR^{10}_{k-e_i,n} \nonumber \\
(\left\langle k+e_j,\omega\right\rangle - \Omega_m)F^{01}_{k+e_j,m} + \bar a_m F^{10}_{k-e_i,n} = -iR^{01}_{k+e_j,m}
\end{eqnarray}
The above three equations have the coefficient matrix of the form
$\left\langle k,\omega\right\rangle I + A_n$ and by the assumption
$(4.2)$ we know that $ |det(\left\langle k,\omega\right\rangle I +
A_n)|\geq \frac{\gamma}{|k|^{\tau}} $ So we get the estimate

\noindent $(1)$: If $n\in \Z^2 _1 \setminus \{{\mathcal L}_1 \cup {\mathcal L}_2\}$, one has
\begin{eqnarray}
|F^{10}_{k,n}|_{\mathcal O}, |F^{01}_{k,n}|_{\mathcal O} \leq c\gamma^{-16}k^{16\tau + 16}max\{
|R^{10}_{k,n}|_{\mathcal O}, |R^{01}_{k,n}|_{\mathcal O} \}
\end{eqnarray}
$(2)$: If $n\in {\mathcal L}_1$ and the corresponding resonant group $£¨m,i,j£©$, one has
\begin{eqnarray}
|F^{10}_{k+e_i,n}|_{\mathcal O}, |F^{10}_{k+e_j,m}|_{\mathcal O} \leq
c\gamma^{-16}k^{16\tau + 16}\max\{ |R^{10}_{k+e_i,n}|_{\mathcal O},
  |R^{10}_{k+e_j,m}|_{\mathcal O} \}
\end{eqnarray}
$(3)$ If $n\in {\mathcal L}_2$, one has
\begin{eqnarray}
|F^{10}_{k-e_i,n}|_{\mathcal O}, |F^{01}_{k+e_j,m}|_{\mathcal O} \leq
c\gamma^{-16}k^{16\tau + 16}\max\{ |R^{10}_{k-e_i,n}|_{\mathcal O},
 |R^{01}_{k+e_j,m}|_{\mathcal O} \}
\end{eqnarray}

\noindent $\bf Solving (4.22)$:
Similarly, we also decompose them into three parts. In this case, the coefficient matrix has the form of
$$
\left\langle k,\omega \right\rangle \pm A_n \otimes I \pm I\otimes A_{n^{\prime}} \quad n,n^{\prime}\in \Z^2 _1
$$
By the assumption:
$$
|\det(\left\langle k,\omega \right\rangle \pm A_n \otimes I \pm
I\otimes A_{n^{\prime}})|\geq \frac{\gamma}{|k|^{\tau}}\quad k\neq 0
$$

\noindent
$(1)$: If $n,n^{\prime}\in \Z^2 _1 \setminus ({\mathcal L}_1 \cup {{\mathcal L}_2})$, one has
\begin{eqnarray}
( \left\langle k,\omega \right\rangle+\Omega_n + \Omega_{n^{\prime}} )F^{20}_{k,n n^{\prime}}
   = -iR^{20}_{k,n n^{\prime}}  \nonumber \\
( \left\langle k,\omega \right\rangle+\Omega_n - \Omega_{n^{\prime}} )F^{11}_{k,n n^{\prime}}
   = -iR^{11}_{k,n n^{\prime}}  \nonumber \\
( \left\langle k,\omega \right\rangle-\Omega_n - \Omega_{n^{\prime}} )F^{02}_{k,n n^{\prime}}
   = -iR^{02}_{k,n n^{\prime}}    \nonumber
\end{eqnarray}
and we get the estimate
\begin{eqnarray}
|F^{20}_{k,n n^{\prime}}|, |F^{11}_{k,n n^{\prime}}|, |F^{02}_{k,n n^{\prime}}| \leq
c\gamma^{-16}k^{16\tau + 16} \max\{ |R^{20}_{k,n n^{\prime}}|, |R^{11}_{k,n n^{\prime}}|, |R^{02}_{k,n n^{\prime}}| \}
\end{eqnarray}
$(2)$: If $n\in \Z^2 _1 \setminus \{{\mathcal L}_1 \cup {\mathcal L}_2\}$, $n^{\prime}\in {\mathcal L}_1$, one has
\begin{eqnarray}
( \left\langle k-e_{i^{\prime}},\omega \right\rangle+\Omega_n - \Omega_{n^{\prime}} )
F^{11}_{k-e_{i^{\prime}},n n^{\prime}} - a_{n^{\prime}}F^{11}_{k-e_{j^{\prime}},n m^{\prime}}
=-iR^{11}_{k-e_{i^{\prime}},n n^{\prime}} \nonumber \\
( \left\langle k-e_{j^{\prime}},\omega \right\rangle+\Omega_n - \Omega_{m^{\prime}} )
F^{11}_{k-e_{j^{\prime}},n m^{\prime}} - a_{m^{\prime}}F^{11}_{k-e_{i^{\prime}},n n^{\prime}}
=-iR^{11}_{k-e_{j^{\prime}},n m^{\prime}} \nonumber
\end{eqnarray}
we have the estimate
\begin{eqnarray}
|F^{11}_{k-e_{i^{\prime}},n n^{\prime}}|_{\mathcal O}, |F^{11}_{k-e_{j^{\prime}},n m^{\prime}}|_{\mathcal O} \leq
c\gamma^{-16}k^{16\tau + 16}\max\{|R^{11}_{k-e_{i^{\prime}},n n^{\prime}}|_{\mathcal O},
|R^{11}_{k-e_{j^{\prime}},n m^{\prime}}|_{\mathcal O}
\}
\end{eqnarray}
\noindent when $n^{\prime} \in {\mathcal L}_2$ is similar, the estimate is similar.

\noindent $(3)$: If $n\in {\mathcal L}_1$,$n^{\prime}\in {\mathcal L}_2$, one has
\begin{eqnarray}
(\left\langle k+e_i + e_{i^{\prime}},\omega \right\rangle+\Omega_n-\Omega_{n^{\prime}})
F^{11}_{k+e_i + e_{i^{\prime}},n n^{\prime}} + a_{n^{\prime}}F^{20}_{k+e_i - e_{j^{\prime}},n m^{\prime}}+
a_{n}F^{11}_{k+e_j + e_{i^{\prime}},m n^{\prime}}
=-i R^{11}_{k+e_i + e_{i^{\prime}},n n^{\prime}} \nonumber \\
(\left\langle k+e_i - e_{j^{\prime}},\omega \right\rangle+\Omega_n+\Omega_{m^{\prime}})
F^{20}_{k+e_i - e_{j^{\prime}},n m^{\prime}} - \bar a_{m^{\prime}}F^{11}_{k+e_i + e_{i^{\prime}},n n^{\prime}}+
a_{n}F^{20}_{k+e_j - e_{j^{\prime}},m m^{\prime}}
=-i R^{20}_{k+e_i - e_{j^{\prime}},n m^{\prime}} \nonumber \\
(\left\langle k+e_j + e_{i^{\prime}},\omega \right\rangle+\Omega_m-\Omega_{n^{\prime}})
F^{11}_{k+e_j + e_{i^{\prime}},m n^{\prime}} + a_{m}F^{11}_{k+e_i + e_{i^{\prime}},n n^{\prime}}+
a_{n}F^{20}_{k+e_i - e_{j^{\prime}},m m^{\prime}}
=-i R^{11}_{k+e_j + e_{i^{\prime}},m n^{\prime}} \nonumber \\
(\left\langle k+e_j - e_{j^{\prime}},\omega \right\rangle+\Omega_m+\Omega_{m^{\prime}})
F^{20}_{k+e_j - e_{j^{\prime}},m m^{\prime}} + a_{m}F^{20}_{k+e_i - e_{j^{\prime}},n m^{\prime}}-
\bar a_{m^{\prime}}F^{11}_{k+e_j + e_{i^{\prime}},m n^{\prime}}
=-i R^{20}_{k+e_j - e_{j^{\prime}},m m^{\prime}} \nonumber
\end{eqnarray}
and we get the estimate
\begin{eqnarray}
&&|F^{11}_{k+e_i + e_{i^{\prime}},n n^{\prime}}|_{\mathcal O},
|F^{20}_{k+e_i - e_{j^{\prime}},n m^{\prime}}|_{\mathcal O},
|F^{11}_{k+e_j + e_{i^{\prime}},m n^{\prime}}|_{\mathcal O},
|F^{20}_{k+e_j - e_{j^{\prime}},m m^{\prime}}|_{\mathcal O}   \\
&&\leq
c\gamma^{-16}|k|^{16\tau + 16}\max
\{
|R^{11}_{k+e_i + e_{i^{\prime}},n n^{\prime}}|_{\mathcal O},
|R^{20}_{k+e_i - e_{j^{\prime}},n m^{\prime}}|_{\mathcal O},
|R^{11}_{k+e_j + e_{i^{\prime}},m n^{\prime}}|_{\mathcal O},
|R^{20}_{k+e_j - e_{j^{\prime}},m m^{\prime}}|_{\mathcal O}
\} \nonumber
\end{eqnarray}
\noindent when $n \in {\mathcal L}_1, n^{\prime} \in {\mathcal L}_1$ or
$n \in {\mathcal L}_2, n^{\prime} \in {\mathcal L}_2$, the estimate is similar.

Now we could give the small-divisor condition in the next step with
new parameters. For simplicity, we only consider the most
complicated case: the second Melnikov condition. Assume that
\begin{eqnarray}
|\det(
\left\langle k,\omega \right\rangle + A_n \otimes I_2 - I_2 \otimes A_{n^{\prime}}
)| > \frac{\gamma}{|k|^{\tau}}  \nonumber
\end{eqnarray}
\noindent then we have
\begin{eqnarray}
&&|\det(
\left\langle k,\omega_+ \right\rangle + A^{+} _n \otimes I_2 - I_2 \otimes A^{+} _{n^{\prime}}
)|   \nonumber \\
&& > |\det(
\left\langle k,\omega \right\rangle + A_n \otimes I_2 - I_2 \otimes A_{n^{\prime}}
)| - c(|k||\hat{\omega}| +
\max \{ |\hat{a}_n|,|\hat{a}_{n^{\prime}}|,|\hat{\Omega}_n|,|\hat{\Omega}_{n^{\prime}}| \}) \nonumber \\
&& > \frac{\gamma}{|k|^{\tau}} - c|k|\varepsilon > \frac{\gamma_{+}}{|k|^{\tau}} \nonumber
\end{eqnarray}
\noindent provided $|k|<K$ where
$K = c(\frac{\gamma-\gamma_+}{\varepsilon})^{\frac{1}{\tau +1}}$. So the small divisor condition in the next step holds automatically for $|k|<K$ and we will deal with other terms in section 6.
\subsection{\textbf{Estimation of coordinate transformation and new perturbation}}
$\\$
\indent With the similar methods in [14], we could get the estimates of
$X_F$ and $\phi^t _F$, just with different parameters.

\begin{Lemma}
Let $D_i = D(r_+ + \frac{i}{4}(r-r_+), \frac{i}{4}s), \quad 0\leq i \leq 4  $. Then we get
\begin{eqnarray}
\|X_F\|_{D_3, \mathcal O} \leq c \gamma^{-16}(r-r_+)^{-c}\varepsilon
\end{eqnarray}
\end{Lemma}

\begin{Lemma}
Let
$\eta = \varepsilon^{\frac{1}{3}},
D_{i\eta}=D(r_+ + \frac{1}{4}(r-r_+), \frac{i}{4}\eta s), \quad
0\leq i \leq 4
$.
If $\varepsilon \ll \gamma^{16} (r-r_+)^c$, then we have
\begin{eqnarray}
\phi^t _F : D_{2\eta}\rightarrow D_{3\eta},\quad -1\leq t \leq 1
\end{eqnarray}
and
\begin{eqnarray}
\|D\phi^t _F - Id  \|_{D_{1\eta}}\leq c\gamma^{-16}(r-r_+)^{-c}\varepsilon
\end{eqnarray}
\end{Lemma}
$\\$
\noindent With above estimates, we could give the estimate of new perturbations. We have
$$
P_+ = \int_0 ^1 \{R(t), F \}\circ \phi^t _F dt + (P-R)\circ \phi^1 _F
$$
where
$
R(t)=R + (1-t)\{N,F\}=(1-t)(N_+ - N)+ tR
$ and
$$
X_{P_+}= \int_0 ^1 (\phi^t _F)^{*}X_{\{R(t),F\}}dt + (\phi^1 _F)^{*}X_{(P-R)}
$$
By Lemma $\bf 4.1$, we get
$$\|D\phi^t _F\|_{D_{1\eta}}\leq 1+ \|D\phi^t _F - Id\|_{D_{1\eta}} \leq 2$$
At the same time, by Cauchy estimate, one has
$$ \|X_{\{R(t),F\}}\|_{D_{2\eta}}\leq c\gamma^{-16}(r-r_+)^{-c}\eta^{-2}\varepsilon^{2}  $$
One the other hand, we have
$$ \|X_{(P-R)}\|_{D_{2\eta}}\leq c\eta \varepsilon $$
To sum up, $P_+$ is bounded by
$$
\|X_{P_+}\|_{D(r_+, s_+)}\leq c\eta \varepsilon + c\gamma^{-16}(r-r_+)^{-c}\eta^{-2}\varepsilon^2 \leq c\varepsilon_+
$$

$\\$
\section{\textbf{Iterative Lemma and Convergence}}
$\\$
\indent For fixed parameters $r,s,\varepsilon, \gamma$, at the $\nu$ th step of the iterative procedure, we define the sequence
\begin{eqnarray}
r_{\nu} &=& r (1-\sum_{i=2} ^{\nu + 1}2^{-i}) \nonumber \\
\varepsilon_{\nu} &=& c \gamma^{-16}(r_{\nu-1} - r_{\nu})^{-c}\varepsilon_{\nu-1}^{\frac{4}{3}} \nonumber \\
\gamma_{\nu} &=& \gamma(1-\sum_{i=2} ^{\nu + 1}2^{-i}) \nonumber \\
\eta_{\nu} &=& \varepsilon_{\nu}^{\frac{1}{3}} \nonumber \\
s_{\nu}&=&\frac{1}{4}\eta_{\nu-1}s_{\nu-1} \nonumber \\
K_{\nu} &=& c(\varepsilon_{\nu} ^{-1} (\gamma_{\nu} - \gamma_{\nu+1}))^{\frac{1}{\tau+1}}
\end{eqnarray}
where $c$ is a constant and the parameters $r_0, s_0, \varepsilon_0, \gamma_0, K_0$ are defined as $r, s, \varepsilon, \gamma, 1$ respectively.

For later use, we define the resonant sets useful for the part of measure estimate:

\begin{eqnarray}
{\mathcal R}^{\nu} =
 \bigcup_{|k| \geq K_{\nu - 1},nm}
 \left(
{\mathcal R}^{\nu} _k \cup
{\mathcal R}^{\nu} _{k,n} \cup
{\mathcal R}^{\nu} _{k,nm}
\right)
\end{eqnarray}

\noindent where each part is defined by

\begin{eqnarray}
{\mathcal R}^{\nu} _k &=&
  \{
  \xi \in {\mathcal O}_{\nu - 1} :
  |\left\langle k,\omega_{\nu} \right\rangle| < \frac{\gamma_{\nu}}{|k|^{\tau}}
  \} \\
{\mathcal R}^{\nu} _{k,n} &=&
  \{
  \xi \in {\mathcal O}_{\nu - 1} :
  |\det(\left\langle k,\omega_{\nu} \right\rangle \pm A^{\nu} _n \otimes I_2 )| < \frac{\gamma_{\nu}}{|k|^{\tau}}
  \}  \\
{\mathcal R}^{\nu} _{k,nm} &=&
  \{
  \xi \in {\mathcal O}_{\nu - 1} :
  |\det(\left\langle k,\omega_{\nu} \right\rangle \pm A^{\nu} _n \otimes I_2 \pm I_2 \otimes A^{\nu} _{m} )| < \frac{\gamma_{\nu}}{|k|^{\tau}}
  \}
\end{eqnarray}

Now we could state the iterative lemma as follows:
\begin{Lemma}
\indent Let $\varepsilon$ is small enough and $\nu \geq 0$, assume that we are at the $\nu$ th step.

$(1)$
$N_{\nu} + \mathcal A_{\nu} + \mathcal B_{\nu} +$ $ {\bar {\mathcal B}}_{\nu}$
is the normal form depending on the parameter $\xi$, where
\begin{eqnarray}
&&N_{\nu} = \left\langle k,\omega_{\nu} \right\rangle + \sum_{n \in \Z^2 _1} \Omega_n^\nu z_n \bar z_n \nonumber \\
&&{\mathcal A}_{\nu} = \sum_{n \in {\mathcal L}_1} a^{\nu} _n e^{i(\theta_i - \theta_j)}z_n \bar z_m \nonumber \\
&&{\mathcal B}_{\nu} = \sum_{n \in {\mathcal L}_2} a^{\nu} _n e^{-i(\theta_i + \theta_j)}z_n z_m \nonumber \\
&&\bar {\mathcal B}_{\nu} = \sum_{n \in {\mathcal L}_2} a^{\nu} _n e^{i(\theta_i + \theta_j)}\bar z_n \bar z_m \nonumber
\end{eqnarray}

\noindent and satisfying the following small divisor conditions:
\begin{eqnarray*}
&&|\left\langle k,\omega_{\nu} \right\rangle| \geq \frac{\gamma_{\nu}}{|k|^{\tau}}  \nonumber \\
&&|\det(\left\langle k,\omega_{\nu} \right\rangle + A^{\nu}_n \otimes I_2)| \geq \frac{\gamma_{\nu}}{|k|^{\tau}} \nonumber \\
&&|\det(\left\langle k,\omega_{\nu} \right\rangle +A^{\nu}_n \otimes
{I_2} \pm I_2 \otimes A^{\nu}_{n^{'}}) | \geq
\frac{\gamma_{\nu}}{|k|^{\tau}}    \nonumber
\end{eqnarray*}
where the matrix is defined as
$$ A^{\nu} _n = \Omega_n \qquad n \in \Z^2 _1 \setminus \left( {\mathcal L}_1 \cup {\mathcal L}_2 \right) $$
$$
A^{\nu} _n =
\left(
\begin{array}{cc}
\Omega^{\nu}_n + \omega^{\nu}_i & a^{\nu}_n \\
a^{\nu}_m & \Omega^{\nu}_m + \omega^{\nu}_j
\end{array}
\right)
\qquad n \in {\mathcal L}_1
$$

$$
A^{\nu} _n =
\left(
\begin{array}{cc}
\Omega^{\nu}_n - \omega^{\nu}_i & -a^{\nu}_n \\
\bar a^{\nu}_m & \Omega^{\nu}_m - \omega^{\nu}_j
\end{array}
\right)
\qquad n \in {\mathcal L}_2
$$

and the parameter $\xi$ is in a closed set $\mathcal O_{\nu}$ of $\R^b$.

$(2)$ $\omega_{\nu},\Omega_{\nu}$ and $a_n^\nu$ are $C_W ^4$ smooth
and satisfy the condition $(\delta=\min\{\bar a-a, \iota\})$
$$
|\omega_{\nu-1} - \omega_{\nu}| < \varepsilon_{\nu-1},
|\Omega_n^{\nu} - \Omega_n^{\nu-1}| <
\varepsilon_{\nu-1}|n|^{-\delta}, |a^{\nu-1}_n - a^{\nu}_n | <
\varepsilon_{\nu-1}|n|^{-\delta}
$$

$(3)$ The perturbation $P_{\nu}$ satisfy condition $(A5)$ and $\|
X_{P_{\nu}}\|_{D(r_{\nu},s_{\nu}),\mathcal O_{\nu}} <
\varepsilon_{\nu}$.

\noindent Then there exists a closed subset ${\mathcal O}_{\nu+1} \subseteq \mathcal O_{\nu}$ defined by
$$
{\mathcal O}_{\nu+1} =
{\mathcal O}_{\nu} \setminus
{\mathcal R^{\nu}}
$$

and a symplectic transformation of variables
$$
\Phi_{\nu}: D(r_{\nu+1},s_{\nu+1}) \times {\mathcal O}_{\nu}
\rightarrow
D(r_{\nu},s_{\nu}) \times {\mathcal O}_{\nu-1}
$$
such that on the  domain $D(r_{\nu+1},s_{\nu+1}) \times {\mathcal
O}_{\nu} $, the Hamiltonian has the form
\begin{eqnarray}
H_{\nu + 1}
&=& N_{\nu+1} + {\mathcal A}_{\nu+1} + {\mathcal B}_{\nu+1} + \bar {\mathcal B}_{\nu+1} \nonumber \\
&=& \left\langle k,\omega_{\nu+1} \right\rangle + \sum_{n \in \Z^2 _1}\Omega_n^{\nu+1} z_n \bar z_n \nonumber \\
&+& \sum_{n \in {\mathcal L}_1}a_n^{\nu+1} e^{i(\theta_i - \theta_j)}z_n \bar z_m \nonumber \\
&+& \sum_{n \in {\mathcal L}_2}a_n^{\nu+1} e^{-i(\theta_i + \theta_j)}z_n  z_m \nonumber \\
&+& \sum_{n \in {\mathcal L}_2}a_n^{\nu+1} e^{i(\theta_i + \theta_j)}\bar z_n \bar z_m \nonumber \\
&+& P_{\nu+1} \nonumber
\end{eqnarray}
with
$$
|\omega_{\nu+1} - \omega_{\nu}| < \varepsilon_{\nu}, \quad
|\Omega^{\nu+1} _n - \Omega^{\nu} _n| < \varepsilon_{\nu}|n|^{-\delta}, \quad
|a^{\nu+1} _n - a^{\nu} _n| < \varepsilon_{\nu}|n|^{-\delta}
$$
The new perturbation $P_{\nu+1}$ satisfy condition $(A5)$ and
$$
\|X_{P_{\nu+1}}\|_{D(r_{\nu+1},s_{\nu+1}), {\mathcal O}_{\nu+1}} < \varepsilon_{\nu+1}
$$
\end{Lemma}

$\\$ Now assume that the assumption of $\bf (A1)-(A5)$ is satisfied.
We could apply the iterative lemma at the $\nu = 0$ step as long as
$\varepsilon_0 , \gamma_0$ are sufficiently small. By an inductive
way, we get the sequence
\begin{eqnarray}
&&{\mathcal O}_{\nu+1} \subseteq {\mathcal O}, \nonumber \\
&&\Phi^{\nu} = \Phi_0 \circ \Phi_1 \circ \cdots \circ \Phi_{\nu} :
D(r_{\nu+1}, s_{\nu+1}) \times {\mathcal O}_{\nu} \rightarrow
D(r_{0}, s_{0}) \times {\mathcal O} \nonumber \\
&&H \circ \Phi^{\nu} = N_{\nu+1} + {\mathcal A}_{\nu+1} + {\mathcal B}_{\nu+1} + \bar {\mathcal B}_{\nu+1} + P_{\nu+1} \nonumber
\end{eqnarray}

Let $\breve{\mathcal O} = \bigcap^{\infty} _{\nu = 0} {\mathcal O}_{\nu}$. It's easy to conclude that
$N_{\nu}, {\mathcal A}_{\nu}, {\mathcal B}_{\nu}, \bar {\mathcal B}_{\nu}, \Phi^{\nu}, D\Phi^{\nu}$ all converge uniformly on $D(\frac{1}{2}r, 0) \times \breve{\mathcal O}$ with
\begin{eqnarray}
&&N_{\infty} = \left\langle \omega_{\infty},I \right\rangle + \sum_{n \in \Z^2 _1}\Omega^{\infty} _n z_n \bar z_n \nonumber \\
&&{\mathcal A}_{\infty} = \sum_{n \in {\mathcal L}_1}a^{\infty} _n e^{i(\theta_i - \theta_j)}z_n \bar z_m \nonumber \\
&&{\mathcal B}_{\infty} = \sum_{n \in {\mathcal L}_2}a^{\infty} _n e^{-i(\theta_i + \theta_j)}z_n z_m \nonumber \\
&&\bar {\mathcal B}_{\infty} = \sum_{n \in {\mathcal L}_2}a^{\infty} _n e^{i(\theta_i + \theta_j)}\bar z_n \bar z_m \nonumber
\end{eqnarray}
\noindent and
\begin{eqnarray}
\varepsilon_{\nu+1}=c\gamma^{-16}(r_{\nu} - r_{\nu+1})^{-c}\varepsilon^{\frac{4}{3}} \rightarrow 0
\end{eqnarray}
provided that $\varepsilon$ is sufficiently small.

Let $\phi^t _H$ be the Hamiltonian flow induced by $X_H$. By $H_{\nu+1} = H \circ \Phi^{\nu}$ one has
$$
\phi^t _H \circ \Phi^{\nu} = \Phi^{\nu} \circ \phi^t _{H_{\nu+1}}
$$
and by the uniform convergence of all the related parameters, we get
$$
\phi^t _H \circ \Phi^{\infty} = \Phi^{\infty} \circ \phi^t _{H_{\infty}}
$$
and
$$
\Phi^{\infty} : D(\frac{1}{2}r,0) \times \breve{\mathcal O} \rightarrow D(r,s) \times \mathcal O
$$
\noindent For parameters $\xi \in \breve{\mathcal O}$, one has
$$
\phi^t _H (\Phi^{\infty}(\T^b \times \{\xi\})) =
\Phi^{\infty} \phi^t _{N_{\infty} + {\mathcal A}_{\infty} + {\mathcal B}_{\infty} + \bar {\mathcal B}_{\infty}}
= \Phi^{\infty} (\T^b \times \{\xi\})
$$
\noindent This means that $\Phi^{\infty} (\T^b \times \{\xi\})$ is an embedded torus which is invariant for the original perturbed Hamiltonian system at $\xi \in \breve{\mathcal O}$. The frequencies $\omega_{\infty}(\xi)$ associated to the tori $\Phi^{\infty} (\T^b \times \{\xi\})$ is slightly different from $\omega(\xi)$. The normal behavior of the invariant torus is governed by normal frequencies $\Omega^{\infty} _n$. \qquad $\Box$

\section{\textbf{Measure Estimate}}
$\\$
\indent Recall the resonant sets at the $\nu$ th step
$
{\mathcal R}^{\nu} =
 \bigcup_{|k| \geq K_{\nu - 1},nm}
 \left(
{\mathcal R}^{\nu} _k \cup
{\mathcal R}^{\nu} _{k,n} \cup
{\mathcal R}^{\nu} _{k,nm}
\right)
$.
To estimate its measure, we need to estimate each single set
$
{\mathcal R}^{\nu} _k, {\mathcal R}^{\nu} _{k,n},  {\mathcal R}^{\nu} _{k,nm}
$
first.

\begin{Lemma}
Fix $|k| \geq K_{\nu - 1},n,m$, one has
$$
meas \left( {\mathcal R}^{\nu} _k \cup {\mathcal R}^{\nu} _{k,n}
\cup {\mathcal R}^{\nu} _{k,nm} \right) < c\frac{\gamma_{\nu}^{\frac
14}}{|k|^{\frac{1}{4}\tau}}
$$
\end{Lemma}
\noindent $\bf Proof$:
One has that
$\omega_{\nu}(\xi) = \omega(\xi) + \sum^{\nu - 1} _{j=0}  P^j _{0l00} (\xi) $
with
$
\sum_{0 \leq j \leq {\nu-1}}|P^j _{0l00}(\xi)|_{{\mathcal O}_{\nu-1}} < \varepsilon
$, and
$
\Omega^{\nu} _n(\xi) = \Omega_n (\xi) +
\sum_{0 \leq j \leq {\nu-1}} P^{011,j}_{nn} (\xi)
$ with
$
\sum_{0 \leq j \leq {\nu-1}} |P^{011,j}_{nn} (\xi)|_{{\mathcal O}_{\nu-1}} < \frac{\varepsilon}{|n|^{\delta}}
$($\delta = \min \{\bar a - a, \iota\} > 0$ ) Similar results also hold for $a_n$. So it's easy to conclude that
$$
\max_{1 \leq j \leq 4} |\partial^j _{\xi} \det \left( \left\langle
k,\omega_{\nu} \right\rangle \pm A_n \otimes I_2 \pm I_2 \otimes A_m
\right) | \geq c|k|
$$
Then the result is obvious. \qquad $\Box$

\begin{Lemma}
The whole measure we need to exclude during the KAM procedure is
\begin{eqnarray}
meas\left( \bigcup_{\nu \geq 0} {\mathcal R}^{\nu}  \right) < c \gamma^{\varsigma}  \qquad \varsigma > 0 \nonumber
\end{eqnarray}
\end{Lemma}

\noindent $\bf Proof$:
Fix one $\nu$ and one $k$, and we only estimate the most complicated term:
\begin{eqnarray}
\bigcup_{n,m}
\left\{
\xi \in {\mathcal O}_{\nu - 1}:
|
\det
\left(
\left\langle k,\omega_{\nu} \right\rangle +
A_n \otimes I_2 -
I_2 \otimes A_m
\right)
| < \frac{\gamma_{\nu}}{|k|^{\tau}}
\right\}  \nonumber
\end{eqnarray}
\noindent Consider its diagonal entry, we only consider one element.
If $|n|^2 - |m|^2 = l \geq c|k|$, then ${\mathcal R}^{\nu}_{k,nm} =
\emptyset$ otherwise we assume $|n| \geq |m|$, by the regularity
property, we get
$$
|\Omega^{\nu} _n - \Omega^{\nu} _m - \varepsilon^{-4} _0 l| < O(|m|^{-\delta}) \qquad \delta=\min\{\bar a - a, \iota\}>0
$$
\noindent so we have that
\begin{eqnarray}
{\mathcal R}^{\nu} _{k,nm} \subseteq {\mathcal Q}^{\nu} _{k,lm} =
\left\{
\xi : |
\det\left(
\left\langle k,\omega_{\nu} \right\rangle +
A_n \otimes I_2 - I_2 \otimes A_m
\right)
| < \frac{\gamma_{\nu}}{|k|^{\tau}} + O(|m|^{-\delta})
\right\} \nonumber
\end{eqnarray}
\noindent it's easy to see that ${\mathcal Q}^{\nu} _{lm} \subseteq {\mathcal Q}^{\nu} _{lm_0}$ for $|m| \geq |m_0|$.
By Lemma $\bf 6.1$ we have
\begin{eqnarray}
meas\left(
\bigcup_{|l|\leq c|k|} \bigcup_{|n|^2 - |m|^2 = l} {\mathcal R}^{\nu} _{k,nm}
\right) &\leq&
\sum_{|l|\leq c|k|}\sum_{|m|\leq |m_0|}meas({\mathcal R}^{\nu} _{k,nm}) +
\sum_{l \leq c|k|} meas({\mathcal Q}_{k,lm_0}) \nonumber \\
&\leq& c \left( \frac{|k|\gamma^{\frac 14}
|m_0|^2}{|k|^{\frac{\tau}{4}}} + O(|k||m_0|^{-\delta}) \right)
\nonumber
\end{eqnarray}
 By choosing appropriate $m_0$ to reach $\frac{\gamma^{\frac 14} |m_0|^2}{|k|^{\frac{\tau}{4}}} = |m_0|^{-\delta}$
(just let $|m_0| = \left( \frac{|k|^{\tau}}{\gamma}
\right)^{\frac{1}{4\delta +8}}$). Then we get the estimate
\begin{eqnarray}
meas \left( \bigcup_{|l|\leq c|k|} \bigcup_{|n|^2 - |m|^2 =
l}{\mathcal R}^{\nu} _{k,nm} \right) <
c\frac{\gamma^{\frac{\delta}{4\delta +
8}}}{|k|^{\frac{\delta\tau}{4\delta + 8}-1}} \nonumber
\end{eqnarray}
\noindent justify the parameter $\tau$ appropriately and we get the result. \qquad $\Box$

\section{\textbf{Appendix}}
$\\$ \indent In this part we give a precise method to construct the
admissible set $ S = \{i_1 = (x_1,y_1),i_2 =(x_2,y_2),\cdots,i_b
=(x_b,y_b)\}   $. It's modified from the appendix in [11] and we
omit some detailed calculation which has been done in [11]. The
points in S will be defined in an inductive way. The first point
$(x_1, y_1) \in \Z^2 _{odd}$ is chosen as $x_1 > b^2, y_1 = 2x_1
^{5^b}$ and the second $x_2 =x_1 ^5  , y_2 =2x_2 ^{5^b} $. If we
have chosen the first $j$ points $i_1, i_2, \cdots,i_j$, then we
define
\begin{eqnarray}
&& x_{j+1} =  x_j ^5 \left(\prod_{1\leq l<m\leq j}((x_m-x_l)^2+(y_m-y_l)^2)+1 \right) \quad 2\leq j \leq b-1                          \nonumber \\
&& y_{j+1} =  2x_{j+1} ^{5^b}    \quad  2\leq j \leq b-1                 \nonumber
\end{eqnarray}
\indent Recall the condition of admissible set (Proposition 1). We verify the conditions one by one.
Given three points $c,d,f \in S$, it's easy to see that
$$
\left\langle c-d,d-f \right\rangle =
(c_1 -d_1)(c_2 - d_2) + (d_1 - f_1)(d_2 - f_2) > 0
$$
So any three points in S can't be three vertices of a rectangle.

To verify condition \textcircled{2}-\textcircled{4}, following the
appendix in [11], it suffices to prove that each equation set in the
following has no integer solution in $\Z^2 _1$ for $c,d,f,g \in S$
and $\{c,d\} \neq \{f,g\}$.

\begin{eqnarray*}
\left\{
\begin{array}
{rcl} &&\left\langle n-g,g-f \right\rangle  =0 \\
&& \left\langle n-c,c-d \right\rangle  =0
\end{array}
\right .
\end{eqnarray*}

\begin{eqnarray*}
\left\{
\begin{array}
{rcl} &&\left\langle n-g,n-f \right\rangle  =0 \\
&& \left\langle n-c,n-d \right\rangle  =0
\end{array}
\right .
\end{eqnarray*}

\begin{eqnarray*}
\left\{
\begin{array}
{rcl} &&\left\langle n-g,g-f \right\rangle  =0 \\
&& \left\langle n-c,n-d \right\rangle  =0
\end{array}
\right .
\end{eqnarray*}

The three equation sets correspond to condition \textcircled{2},\textcircled{3},\textcircled{4} respectively and we denote them by $\bf \Rmnum{1},\Rmnum{2},\Rmnum{3}$ respectively.

For $\bf \Rmnum{1}$,
$\\$
\indent $\bf \Rmnum{1}-\rmnum{1} :$ If only one element of $\{|c|,|d|,|f|,|f|\}$ reaches the maximun value of them.

$\bf \Rmnum{1}-\rmnum{1}-(1) :$ $|d|$ or $|f|$ reaches the maximum. Without lose of generality, just assume $d$.
It's easy to get that
\begin{eqnarray}
n_2 = c_2 + \frac
{
(g_1-c_1)(g_1-f_1)(c_1-d_1)+(g_2-c_2)(c_1-d_1)(g_2-f_2)
}
{
(c_1-d_1)(g_2-f_2)-(c_2-d_2)(g_1-f_1)
} \nonumber
\end{eqnarray}
By the fact $d_2 \gg$ other ones, the second term in the right side couldn't be an integer, implying $n_2 \notin \Z$, which is a contradiction.

$\bf \Rmnum{1}-\rmnum{1}-(2):$ $|g|$ or $|c|$ reaches the maximum. Without lose of generality, just assume $|g|$.
By calculation, we get
\begin{eqnarray}
n_2 &=&
g_2 +
2(
c_1 ^{5^b-1}+c_1 ^{5^b-2}d_1+\cdots+c_1d_1 ^{5^b-2}+d_1 ^{5^b-1}
)(g_1-f_1)  \nonumber \\
&+&\frac
{
g_1-c_1+4(g_1-f_1)A^2 +
4(f_1 ^{5^b}-c_1 ^{5^b})A
}
{
2(g_1 ^{5^b-1}+g_1 ^{5^b-2}f_1+\cdots+g_1f_1 ^{5^b-2}+f_1 ^{5^b-1})- 2A
}\nonumber
\end{eqnarray}
where
$A= (c_1 ^{5^b-1}+c_1 ^{5^b-2} d_1+ \cdots +c_1d_1 ^{5^b - 2}+d_1 ^{5^b-1})$
By the fact $g_1 \gg$ other ones, we know that the last term is not an integer, so $n_2 \notin \Z$, which is a contradiction.

$\bf \Rmnum{1}-\rmnum{2}:$ If two elements of $|c|,|d|,|f|,|g|$ reach the maximum, by the structure of S we know that these two points must be equal, and the case when three points reach the maximum will not happen.

$\bf \Rmnum{1}-\rmnum{2}-(1):$ when $|d|=|g|$ reach the maximum, we have $d=g$.
By calculation we get
\begin{eqnarray}
n_2=
\frac
{
4g_1 ^{5^b}A - 4c_1 ^{5^b}B + (g_1-c_1)
}
{
4(A-B)
} \nonumber
\end{eqnarray}
\noindent where $A= (g_1 ^{5^b-1} + g_1 ^{5^b-2}f_1 +\cdots+g_1f_1 ^{5^b-2}+f_1 ^{5^b-1}) $
and $B= (g_1 ^{5^b-1} + g_1 ^{5^b-2}c_1 +\cdots+g_1c_1 ^{5^b-2}+c_1 ^{5^b-1}) $
without losing generality, we assume $|c|<|f|$(the case $|c|=|f|$ will not happen).
According to the structure of S, it's easy to see that in the expression above, the denominator is divisible by $|c_1|^4$, in the numerator, all terms are divisible by $|c_1|^4$ except for $c_1$. It means that $n_2 \notin \Z$ which is a contradiction.

$\bf \Rmnum{1}-\rmnum{2}-(2):$ If $|d|=|f|$ reach the maximum, we have $d=f$, then we have
$$
\left\langle n-g,g-f \right\rangle = \left\langle n-g,g-d \right\rangle = \left\langle n-c,c-f \right\rangle = 0
$$
In this case $c,d,g$ are three vertices of a rectangular, which is a contradiction.

$\bf \Rmnum{1}-\rmnum{2}-(3):$ If $|c|=|g|$ reach the maximum, we have $c=g$. From
$$
\left\langle n-g,g-f \right\rangle = \left\langle n-g,g-d \right\rangle = \left\langle n-c,c-f \right\rangle = 0
$$
we know $d,f,g$ lie on the same line, which is a contradiction.

Now we turn to $\bf \Rmnum{3}$.

$\bf \Rmnum{3}-\rmnum 1:$ If only one element of $|c|,|d|,|f|,|g|$ reaches the maximum value of them and each one is different from others.

$\bf \Rmnum{3}-\rmnum 1-(1):$ $|d|$ reaches the maximum.

We have
$$
\left\langle n-c,n-d \right\rangle = 0 \qquad
\left\langle n-g,g-f \right\rangle = 0
$$
We take $g$ as the origin and from above we get an equation about $n_1$.
\begin{eqnarray}
(f_1 ^2+f_2 ^2)n_1 ^2+\left((c_2+d_2)f_1f_2-(c_1+d_1)f_2 ^2 \right)n_1+f_2 ^2(c_1d_1+c_2d_2)=0 \nonumber
\end{eqnarray}
The discriminant of the equation is
\begin{eqnarray}
\Delta=
\left(
f_1f_2d_2+
\left(
c_2f_1f_2-c_1f_2 ^2-d_1f_2 ^2-\frac{2c_2f_2(f_1 ^2+f_2 ^2)}{f_1}
\right)-\alpha
\right)^2 \nonumber
\end{eqnarray}
where $\alpha \sim \frac{d_1}{d_2} \ll \frac{1}{|f_1|}$ so we get
$$
n_1=\frac{-\left((c_2+d_2)f_1f_2-(c_1+d_1)f_2 ^2\right)\pm \sqrt{\Delta}}
{2(f_1 ^2+f_2 ^2)}
$$
while $\sqrt{\Delta}=$ some integer plus $\alpha$, so the numerator in the above expression is not an integer. So we conclude that $n_1 \notin \Z$, which is a contradiction.

$\bf \Rmnum{3}-\rmnum 1-(2)$ $|f|$ reaches the maximum.

As before, we take $g$ as the origin and we get
\begin{eqnarray}
(f_1 ^2+f_2 ^2)n_1 ^2+
\left(
(c_2+d_2)f_1f_2-(c_1+d_1)f_2 ^2
\right)n_1 +
(c_1d_1+c_2d_2)f_2 ^2 = 0 \nonumber
\end{eqnarray}
The discriminant is
\begin{eqnarray}
\Delta=(c_1-d_1)^2f_2 ^4-4c_2d_2f_2 ^4-2(c_1+d_1)(c_2+d_2)f_1f_2 ^3
+(c_2-d_2)^2f_1 ^2f_2 ^2-4c_1d_1f_1 ^2f_2 ^2 <0 \nonumber
\end{eqnarray}
by the fact $c_2d_2 \gg (c_1-d_1)^2$ and $f_2 \gg$ other terms.

$\bf \Rmnum{3}-\rmnum 1-(3)$ $|g|$ reaches the maximum.

It's easy to see $|n|^2 > |g|^2-|f|^2$ and we get
$$
|n|^2 + |m|^2 -|c|^2-|d|^2 > 0
$$

$\bf \Rmnum{3}-\rmnum{2}$ If only one element of $\{|c|,|d|,|f|,|g|\}$ reaches the maximun and two of the others are the same.

$\bf \Rmnum{3}-\rmnum{2}-(1)$ $|d|$ reaches the maximum, and $c=g$, without losing generality, we just assume $c=g$. We could assume $|f|>|c|$, the case $|f|<|c|$ is similar.

In this case, by calculation we get
\begin{eqnarray}
n_2=
c_2+
\frac
{
(f_1-c_1)^2d_2-(f_1-c_1)(f_2-c_2)d_1
}
{
(f_1-c_1)^2+(f_2-c_2)^2
}
+
\frac
{
2c_1f_1(f_1 ^{5^b-2}+f_1 ^{5^b-3}c_1+\cdots+f_1c_1 ^{5^b-3}+c_1 ^{5^b-2})
}
{
1+4(f_1 ^{5^b-1}+f_1 ^{5^b-2}c_1+\cdots+f_1c_1 ^{5^b-2}+c_1 ^{5^b-1})^2
} \nonumber
\end{eqnarray}
By the structure of S, we could conclude that the term
$$
\frac{(f_1-c_1)^2d_2-(f_2-c_2)(f_1-c_1)d_1}
{(f_1-c_1)^2+(f_2-c_2)^2}
$$
must has a form of an integer plus
$
\frac{2(f_1-c_1)^2f_1 ^{\kappa}-(f_2-c_2)(f_1-c_1)f_1 ^{\tilde{\kappa}}}
{(f_1-c_1)^2 + (f_2-c_2)^2}
$ where $0<\kappa\leq 5^b$ and $5\leq \tilde{\kappa} \leq 5^{b-1}$.
Then we have
$
\frac{(f_1 - c_1)^2 f^{\kappa} _1}{(f_1-c_1)^2+(f_2-c_2)^2} \leq \frac{1}{f_1 ^{5^b-2}}
$
and
$
\frac{c_1f_1(f_1 ^{5^b-2}+f_1 ^{5^b-3}c_1+\cdots+f_1c_1 ^{5^b-3}+c_1 ^{5^b-2})}
{1+(f_1 ^{5^b-1}+f_1 ^{5^b-2}c_1+\cdots+f_1c_1 ^{5^b-2}+c_1 ^{5^b-1})^2}
\leq \frac{1}{f_1 ^{5^b-2}}
$ while
$\frac{(f_2-c_2)(f_1-c_1)f_1 ^{\tilde{\kappa}}}
{(f_2-c_2)^2+(f_1-c_1)^2} > f_1 ^{6-5^b}
$ which is much lager than the former two, and these three terms all $\ll 1$,
hence $n_2 \notin \Z$.
The case $c=f$ is similar to $\bf \Rmnum{3}-\rmnum{1}-(1)$.

$\bf \Rmnum{3}-\rmnum{2}-(2)$ $|f|$ reaches the maximum and $g=c$ or $g=d$. Without losing generality, we just assume $g=c$.

We have
$$
\left\langle n-c,n-d \right\rangle=0 \qquad \left\langle n-c,c-f \right\rangle=0
$$
If we take $c$ as the origin, we get
$
n_2=\frac{f_1 ^2d_2-f_1f_2d_1}{f_1 ^2+f_2 ^2} \notin \Z
$

$\bf \Rmnum{3}-\rmnum{2}-(3)$ $|g|$ reaches the maximum. It's similar to $\bf \Rmnum{3}-\rmnum{1}-(3)$.

$\bf \Rmnum{3}-\rmnum{3}$ If two elements of $\{|c|,|d|,|f|,|g|\}$ reach the maximum of their values.

$\bf \Rmnum{3}-\rmnum{3}-(1)$ $|g|=|d|$ reach the maximum. We have $g=d$. Then we get
$$
n_1 = c_1+d_1-f_1+
\frac
{
(d_1-f_1)^2(d_1-c_1)+(f_2-c_2)(d_1-f_1)(d_2-f_2)-(d_1-f_1)^3
}
{
(d_1-f_1)^2+(d_2-f_2)^2
}
$$
which implies that $n_1 \notin \Z$ due to the fact $d_2 \gg$ other terms.

$\bf \Rmnum{3}-\rmnum{3}-(2)$ $|f|=|d|$ reach the maximum. We have $f=d$ and take $g$ as the origin as before. Then we get
$$
(d_1 ^2+d_2 ^2)n_1 ^2+(d_1d_2c_2-c_1d_2 ^2)n_1+c_1d_1d_2 ^2+c_2d_2 ^3 = 0
$$
The discriminant is
$$
\Delta =
(d_1d_2c_2-c_1d_2 ^2)^2-4(d_1 ^2+d_2 ^2)(c_1d_1d_2 ^2+c_2d_2 ^3)<0
$$
due to the fact $d_2 \gg$ other terms. So $n_1$ doesn't exist, which is a contradiction.

At last we turn to $\bf \Rmnum{2}$.

We know $\{c,d\}\notin \{f,g\}$. If $\{c,d\}\cap \{f,g\} \neq \emptyset$, without losing generality we just assume $c=f$. The equation becomes
$$
\left\langle n-g,g-d \right\rangle=0
\left\langle n-c,n-d \right\rangle=0
$$
We have proved that it has no solution in $\Z^2$ in $\bf \Rmnum{3}-(\rmnum{3})-(2)$.

Now we concentrate on the remaining case when the four elements are different from each other. Without losing generality, we just assume $|d|$ reach the maximum of their values. It's easy to see $|n|^2 \ll d_1$ and
$$
n_2=c_2+
\frac{(c_1-n_1)d_1-c_2(c_2-f_2-g_2)-\left\langle f,g \right\rangle-n_1(c_1-f_1-g_1)}
{c_2+d_2-f_2-g_2}
$$
If $c_1=n_1$, then in the above expression the numerator is smaller than $d_1$, and if $c_1 \neq n_1$, the numerator is smaller than $\frac{d_2}{2}$, we still have $n_2 \notin \Z$, which is a contradiction.

\thispagestyle{empty}

\end{document}